\documentclass[11pt,oneside,reqno]{amsart}
\usepackage[utf8]{inputenc}
\usepackage[T1]{fontenc}
\usepackage{lmodern}
\usepackage{amssymb, amsmath, amsthm, mathtools, cases, enumerate}
\usepackage{color}
\usepackage{hyperref,mathscinet}
\usepackage{natbib}
\usepackage{microtype}
\usepackage{hyphenat}
\usepackage{enumitem}
\usepackage{amsfonts}       

\hypersetup{hidelinks}

\hypersetup{
	unicode=true,
	breaklinks=true,
	pdftitle={Optimal function spaces in weighted Sobolev embeddings with $\alpha$-homogeneous weights},
	pdfauthor={Ladislav Drážný}
	}

\title[Optimal spaces in embeddings with $\alpha$-homogeneous weights]{Optimal function spaces in weighted Sobolev embeddings with $\alpha$-homogeneous weights}
\author{Ladislav Drážný}

\address{Ladislav Drážný, Charles University, Faculty of Mathematics and Physics, Department of Mathematical Analysis, Sokolovská~49/83, 186~75 Praha~8, Czech Republic; Czech Technical University in Prague, Faculty of Electrical Engineering, Department of Mathematics, Technická~2, 166~27 Praha~6, Czech Republic}
\email{drazny@karlin.mff.cuni.cz}
\urladdr{\href{https://orcid.org/0009-0003-1807-6003}{0009-0003-1807-6003}}

\date{\today}

\newenvironment{myproof}{
  \par\medskip\noindent
  \textit{Proof}.
}{
\newline
\rightline{$\qedsymbol$}
}

\newenvironment{myproof3}{
  \par\medskip\noindent
  \textit{Proof of Proposition \ref{d2}}.
}{
\newline
\rightline{$\qedsymbol$}
}

\newenvironment{myproof4}{
  \par\medskip\noindent
  \textit{Proof of Proposition \ref{PS}}.
}{
\newline
\rightline{$\qedsymbol$}
}

\theoremstyle{plain}
\newtheorem{thm}{Theorem}
\numberwithin{thm}{section}
\newtheorem{lemma}[thm]{Lemma}

\newtheorem{prop}[thm]{Proposition}

\theoremstyle{remark}

\newtheorem{rem}[thm]{Remark}

\renewcommand\qedsymbol{\square}

\newcommand{\R}{\mathbb{R}}
\newcommand{\N}{\mathbb{N}}

\newcommand{\ballm}{B_{\mu}}
\newcommand{\rn}{\mathbb{R}^n}
\newcommand{\RNum}[1]{\uppercase\expandafter{\romannumeral #1\relax}}

\numberwithin{equation}{section}

\thanks{This research was partly supported by grant No.~23-04720S of the Czech Science Foundation and by the Primus research programme PRIMUS/21/SCI/002 of Charles University.}

\begin{document}
\setcitestyle{numbers}
\bibliographystyle{plainnat}

\begin{abstract}
We study weighted Sobolev inequalities on open convex cones endowed with $\alpha$-homogeneous weights satisfying a certain concavity condition. We establish a so-called reduction principle for these inequalities and characterize optimal re\-ar\-range\-ment\--in\-vari\-ant function spaces for these weighted Sobolev inequalities. Both optimal target and optimal domain spaces are characterized. Abstract results are accompanied by general yet concrete examples of optimal function spaces. For these examples, the class of so-called Lorentz--Karamata spaces, which contains in particular Lebesgue spaces, Lorentz spaces, and some Orlicz spaces, is used.
\end{abstract}
\maketitle


\section{Introduction}\label{beg}
Various Sobolev spaces and inequalities have played important roles in mathematics for decades. Their applications include (but are not limited to) analysis of partial differential equations, calculus of variations, or harmonic analysis. In this paper, we study weighted Sobolev inequalities on open convex cones $\Sigma$ in $\R^n$ with $\alpha$-homogeneous weights $w \colon \overline{\Sigma} \rightarrow [0, \infty)$ and establish their optimal, in a sense, versions. Throughout the entire paper, we assume that $n\in\N$, the dimension of $\R^n$, is greater than $1$. The weight $w$ is a nonnegative continuous $\alpha$-homogeneous function for some $\alpha > 0$ such that $w^{\frac{1}{\alpha}}$ is concave in $\Sigma$. We define the weighted measure $\mu$ on $\Sigma$ as
\begin{equation*}
\mu(E) = \int\limits_{E} w(x) \, dx
\end{equation*}
for every Lebesgue measurable set $E \subseteq \Sigma$.

An important example of an admissible weight is the monomial weight defined as $w(x)=x_1^{A_1} \cdots x_k^{A_k}$ on the cone $\Sigma = \{x_i>0, i = 1, \dots, k\}$ for some  $1\leq k \leq n$ and $A_i > 0$, which is $\alpha=A_1+\cdots+A_k$ homogeneous. These monomial weights and corresponding weighted Sobolev inequalities have been quite fashionable recently (e.g., \cite{C:17, C:21, GH:21, L:17, W:23}, to name a few). Their importance was observed in particular in \cite{CR-O:13b}, where they were used in connection with the regularity of stable solutions to certain planar reaction--diffusion problems, and they were studied in more detail in \cite{CRO}. In particular, in the latter, they established a new isoperimetric inequality for the monomial weights. What is quite surprising and remarkable is that the isoperimetric quotient is minimized by intersections of balls with the cone $\Sigma$ despite the fact that the monomial weights are not radially symmetric. Using their isoperimetric inequality, they also established the corresponding Sobolev inequality, namely for $p\in[1, D)$
\begin{equation*}
\|u\|_{L^{p^*}(\Sigma, \mu)} \leq C_{D,p} \|\nabla u\|_{L^p(\Sigma, \mu)} \quad \text{for every $u\in \mathcal C_c^1(\R^n)$},
\end{equation*}
where $D = n + k$ and $p^* = (Dp)/(n-p)$. Notice that the usual role of the dimension $n$ is replaced by the sum of the dimension and the order of homogeneity of the weight. Furthermore, note that the presence of $\R^n$ (instead of $\Sigma$) in $C_c^1(\R^n)$ is not a typo, because the functions are not required to vanish on the boundary of $\Sigma$.

Later, the isoperimetric inequality for the monomial weights was greatly generalized in~\cite{CR-OS} (see also~\cite{CGPR-O:22}). In particular, their general isoperimetric inequality applies to the general setting of convex cones with $\alpha$-homogeneous weights considered in this paper. Noteworthily, arguments based on symmetrization can often be successfully used even though the weights need not be radially symmetric. Besides the monomial weights, some examples of admissible weights and cones are the following. If $ \Sigma = (0, \infty)^{n}$, then we can consider the weight $w(x) = (A_{1}x_{1}^{\frac{1}{r}} + \dots + A_{n}x_{n}^{\frac{1}{r}})^{\alpha r}$ for $r \geq 1, A_{i} \geq 0$, and $\alpha > 0$. For $\Sigma\neq \R^n$, the weight $w(x) = \operatorname{dist}(x, \partial \Sigma)^{\alpha}$, $\alpha > 0$ is admissible. More examples of admissible weights can be found in \cite{CR-OS}.  It is worth pointing out that it may happen that the weight $w$ does not vanish on $\partial \Sigma$. Such weights can be obtained for instance by restricting an admissible weight $w$ to a subcone of $\Sigma$.

In this paper, we study Sobolev inequalities of the form
\begin{align}\label{intro:Sob_weighted_ri}
\|u\|_{Y(\Sigma, \mu)} \leq C \| \nabla^m u \|_{X(\Sigma, \mu)} \quad \text{for every $u\in V_0^m X(\Sigma, \mu)$},
\end{align}
where (see~Section~\ref{intro} for precise definitions) $X(\Sigma, \mu)$ and $Y(\Sigma, \mu)$ are re\-ar\-range\-ment\--in\-vari\-ant function spaces on $\Sigma$ endowed with the measure $\mu$ and $V_0^m X(\Sigma, \mu)$ is a suitable $m$th order ($1\leq m < D$) Sobolev-type space built upon $X(\Sigma, \mu)$. We study the question of optimal function spaces in \eqref{intro:Sob_weighted_ri}. On the one hand, for a given re\-ar\-range\-ment\--in\-vari\-ant function space $X(\Sigma, \mu)$, we will describe the optimal re\-ar\-range\-ment\--in\-vari\-ant function space $Y(\Sigma, \mu)$ with which the inequality (\ref{intro:Sob_weighted_ri}) holds (see Theorem~\ref{zopt}). The optimality is understood in the sense that $Y(\Sigma, \mu)$ cannot be replaced by a strictly smaller re\-ar\-range\-ment\--in\-vari\-ant function space. On the other hand, for a given $Y(\Sigma, \mu)$, we will describe the optimal re\-ar\-range\-ment\--in\-vari\-ant function space $X(\Sigma, \mu)$ with which the inequality (\ref{intro:Sob_weighted_ri}) is valid (see Theorem~\ref{zopt22}). This time, the optimality is understood in the sense that $X(\Sigma, \mu)$ cannot be replaced by a strictly larger re\-ar\-range\-ment\--in\-vari\-ant function space.

The point of departure for us is the suitable isoperimetric inequality obtained in \cite{CR-OS}, which we exploit to obtain a suitable weighted Pólya--Szeg\H{o} inequality (see Proposition~\ref{PS}). By utilizing it and combining it (when $m>1$) with the clever idea of iteration (see~\cite{CPS}), we prove a so-called reduction principle for the inequality \eqref{intro:Sob_weighted_ri} (see~Theorem~\ref{RP}). To establish these results, we need to exploit a lot of different techniques developed and improved over time together with results from both classical and contemporary theory of (re\-ar\-range\-ment\--in\-vari\-ant) function spaces. The question of optimal re\-ar\-range\-ment\--in\-vari\-ant function spaces for a large number of various Sobolev inequalities has been intensively studied for more than two decades (e.g., \cite{ACPS:18, CP2, CCP, CPS, CPS:20, EKP:00, KP:06, M, M:21b}). Nevertheless, the setting considered in this paper appears not be covered, and so we aim to fill this gap in the actively developing setting of Sobolev inequalities on convex cones with $\alpha$-homogeneous weights. Furthermore, abstract theorems are accompanied by several concrete (yet substantially general) examples of optimal function space in the inequality \eqref{intro:Sob_weighted_ri} (see~Section~\ref{ex}).

The class of re\-ar\-range\-ment\--in\-vari\-ant function spaces is very general and contains a large number of classical function spaces. We will briefly informally introduce this class of function spaces here. Loosely speaking, re\-ar\-range\-ment\--in\-vari\-ant function spaces are often suitable for measuring integrability and, arguably, their most fundamental feature is that re\-ar\-range\-ment\--in\-vari\-ant function norms depend only on the measure of level sets. By that we mean that if $u$ and $v$ are two measurable functions such that the measures of the sets $\{x\colon |u(x)| > \lambda\}$ and $\{x\colon |v(x)| > \lambda\}$ are the same for every $\lambda > 0$, then their norms are equal. For example, the rearrangement invariance of the $L^p$ norms follows from the well-known layer cake representation formula (\cite[Theorem~1.13]{LiebLoss2001}). Apart from the Lebesgue spaces, other well-known examples of re\-ar\-range\-ment\--in\-vari\-ant function spaces are Lorentz spaces $L^{p,q}$ or Orlicz spaces $L^A$. Orlicz spaces naturally generalize Lebesgue spaces and their usefulness stems from their usage in analysis of nonlinear partial differential equations and variational problems whose nonlinearity is non-polynomial. Very loosely speaking, they measure integrability in a more fine-grained way by replacing power functions by general convex functions (i.e., $|f(x)|^p$ versus $A(|f(x)|)$), which allows for capturing non-polynomial growth. Lorentz spaces, which contain the famous weak Lebesgue spaces (corresponding to $q=\infty$), are not only intimately connected with the interpolation theory and harmonic analysis (e.g., \cite{BL:76, G1, SW:71}) but also with the theory of Sobolev spaces. For example, assuming $n>1$, while a function whose weak gradient (locally) belongs to the Lebesgue space $L^n(\R^n)$ need not be (locally) bounded or be differentiable (in the classical sense) at any point, the situation changes dramatically when its weak gradient (locally) belongs to the Lorentz space $L^{n,1}(\R^n)$, which is a slightly (yet essentially) smaller function space than $L^{n}(\R^n)$. In the latter case, not only is the function (locally) bounded, but it has a continuous representative that is differentiable  almost everywhere (see~\cite{S:81}). Another well-known result is that the classical Sobolev embedding into the critical Lebesgue space $L^{p^*}$ can be improved by replacing $L^{p^*}$ with the Lorentz space $L^{p^*, p}\subsetneq L^{p^*}$ (e.g., \cite{O:63, P:66, T:98}). This improvement can also be essential in analysis of partial differential equations and variational problems with critical growth (see~\cite{R:06}).

A considerably general (yet still reasonably concrete) subclass of re\-ar\-range\-ment\--in\-vari\-ant function spaces is constituted by so-called Lorentz--Karamata spaces (see~\cite{MR2091115, EKP:00, MR1927106, P2}). Among other more delicate function spaces, it contains not only Lebesgue spaces and Lorentz spaces but also a lot of Orlicz spaces (in particular, those of ``exponential'' and ``logarithmic'' type (see also~\cite{OP:99})). We will consider this class of function spaces in Section~\ref{ex}, where we provide concrete examples of optimal function norms in \eqref{intro:Sob_weighted_ri}. With this choice, not only can we describe the optimal function spaces explicitly, but it also enables us to capture delicate integrability properties in limiting situations (in particular, loosely speaking when $X$ is ``close to'' $L^{\frac{D}{m}}$ (cf.~\cite{BW:80, H:79, M:71, T:67})).

\section{Preliminaries} \label{intro}
In the whole paper we use the convention $\frac{1}{\infty} = 0$ and $0\cdot \infty = 0$. When $E\subseteq (0, \infty)$ is Lebesgue measurable, we denote by $\lambda(E)$ its Lebesgue measure.

Let $(R, \mu)$ be a $\sigma$-finite nonatomic measure space. By $\mathcal{M}(R, \mu)$ we will denote the class of all $\mu$-measurable functions on $R$ whose values lie in $\mathbb{R} \cup \{- \infty, \infty\}$. We will denote the class of all $\mu$-measurable functions on $R$ whose values lie in $[0, \infty]$ by $\mathcal{M}^{+}(R, \mu)$. And the class of all functions in $\mathcal{M}(R, \mu)$ that are finite $\mu$-almost everywhere in $R$ will be denoted by $\mathcal{M}_{0}(R, \mu)$. 

Now we introduce re\-ar\-range\-ment\--in\-vari\-ant Banach function spaces and their basic poperties. The theory that is presented here except of the Sobolev spaces follows the first three chapters of \cite{BS}.

Let $f \in \mathcal{M}(R, \mu)$. The \emph{nonincreasing rearrangement} of the function $f$ is the function $f^{\ast}_{\mu} \colon (0, \infty) \to [0, \infty]$ defined by
\begin{align*}
f^{\ast}_{\mu}(t) = \inf\limits_{\lambda > 0}\left(\mu(\{x \in R; \lvert f(x) \rvert > \lambda \}) \leq t \right) \, \, \, \, t \in (0, \infty).
\end{align*}
The function $f^{\ast}_{\mu}$ is nonincreasing and right-continuous. If $f \in \mathcal{M}(R, \mu)$, $g \in \mathcal{M}(S, \nu)$ satisfy
\begin{align}\label{rpp5}
f^{\ast}_{\mu}(t) = g^{\ast}_{\nu}(t), \, \, \, \, t \in (0, \infty),
\end{align}
we say that $f$ and $g$ are \emph{equimeasurable}. For instance, the functions $f$ and $f^{\ast}_{\mu}$ are equimeasurable.

A mapping $\rho$ on $\mathcal{M}^{+}(R, \mu)$ with values in $[0, \infty]$ we call a \emph{ re\-ar\-range\-ment\--in\-vari\-ant Banach function norm} if all the following properties are satisfied for all $f, g \in \mathcal{M}^{+}(R, \mu)$, $\{f_{k}; k \in \mathbb{N}\} \subseteq \mathcal{M}^{+}(R, \mu)$, $c \in [0, \infty)$ and $A \subseteq R$ such that A is $\mu$-measurable.
\begin{enumerate}
  \item\label{th1} the \emph{norm axiom}: $\rho(f) = 0$ if~and~only~if $f = 0$ $\mu$-almost everywhere in $R$, $\rho(cf) = c\rho(f)$, $\rho(f + g) \leq \rho(f) + \rho(g)$;
  \item\label{th2} the \emph{lattice axiom}: if $g \leq f$ $\mu$-almost everywhere in $R$, then $\rho(g) \leq \rho(f)$;
  \item\label{th3} the \emph{Fatou axiom}: if $f_{k} \uparrow f$ $\mu$-almost everywhere in $R$, then $\rho(f_{k}) \uparrow \rho(f)$;
  \item\label{th4} the \emph{nontriviality axiom}: if $\mu({A}) < \infty$, then $\rho(\chi_{A}) < \infty$;
  \item\label{th5} the \emph{local embedding in $L^{1}$}: if $\mu(A) < \infty$, then
\begin{align}\label{emb}
\int_{A} f \, d\mu \leq K_{A}\rho(f),
\end{align}
where $K_{A} \geq 0$ is a constant which may depend on $A$ but which does not depend on $f$;
\item the \emph{re\-ar\-range\-ment\--in\-vari\-ance axiom}:  if $f^{\ast}_{\mu} = g^{\ast}_{\mu}$, then $\rho(f) = \rho(g)$.
\end{enumerate}
The collection of all functions $f \in \mathcal{M}(R, \mu)$ such that $\rho(\lvert f \rvert) < \infty$ is called a \emph{re\-ar\-range\-ment\--in\-vari\-ant Banach function space}. We will denote it by $X(\rho)$, $X(R, \mu)$ or just by $X$.

As their name suggests, re\-ar\-range\-ment\--in\-vari\-ant Banach function spaces are Banach spaces. Textbook examples of re\-ar\-range\-ment\--in\-vari\-ant spaces are the Lebesgue spaces $L^p(R, \mu)$ for $p\in[1, \infty]$. Every re\-ar\-range\-ment\--in\-vari\-ant space contains simple functions (i.e., linear combinations of characteristic functions of $\mu$-measurable sets of finite measure) and is contained in $\mathcal{M}_{0}(R, \mu)$.

Now we define suitable weighted Sobolev spaces built on re\-ar\-range\-ment\--in\-vari\-ant spaces. We start with some notation and conventions used in the rest of this paper. Throughout the rest of this paper, we assume that $n\in\N$, $n\geq2$, is the dimension of $\R^n$. Furthermore, we assume that $\Sigma$ is a nonempty open convex cone with vertex at the origin, i.e., $\Sigma \subseteq \R^{n}$ is a nonempty open convex set such that for every $x \in \Sigma$ and for every $r > 0$, we have $rx \in \Sigma$. We also assume that $w \colon \overline{\Sigma} \rightarrow [0, \infty)$ is a nonnegative continuous function that is not identically zero. Furthermore, we assume that $w$ is $\alpha$-homogeneous for some $\alpha > 0$, i.e., for every $x \in \overline{\Sigma}$ and for every $s > 0$, we have $w(s x) = s^{\alpha} w(x)$. Finally, we assume that the function $w^{\frac{1}{\alpha}}$ is concave in $\Sigma$. We set 
\begin{align}\label{defD}
D = n + \alpha. 
\end{align}
We assume that $m\in\N$ is such that 
\begin{align}\label{defm}
1\leq m < D.
\end{align}

We define the weighted measure $\mu$ on $\Sigma$ as
\begin{equation*}
\mu(E) = \int\limits_{E} w(x) \, dx
\end{equation*}
for every Lebesgue measurable set $E \subseteq \Sigma$.

We now introduce the Sobolev spaces that we will work with. Let $k \in \mathbb{N}$ and let $u \colon \Sigma \rightarrow \R $ be $k$-times weakly differentiable function in $\Sigma$ (i.e., it has all weak derivatives up to the $k$-th order). We denote by $\nabla^l u$, $l \in \{1, \dots, k\}$, the vector of all $l$-th order weak derivatives of $u$. We also set $\nabla^0 u = u$.

Let $X(\Sigma, \mu)$ be a re\-ar\-range\-ment\--in\-vari\-ant space. We say that $u$ belongs to the space $V^{k}X(\Sigma, \mu)$ if
\begin{equation*}
\left\lvert \nabla^{k} u \right\rvert \in X(\Sigma, \mu).
\end{equation*}
We say that $u$ belongs to the space $V^{k}_{0}X(\Sigma, \mu)$ if $u \in V^{k}X(\Sigma, \mu)$ and for every $l \in \{0, 1, \dots, k-1 \}$ and for every $\lambda > 0$ it holds that
\begin{equation*}
\mu\left(\left\{x \in \Sigma; \left\lvert \nabla^{l} u(x) \right\rvert > \lambda \right\}\right) < \infty.
\end{equation*}
For short, we will write $\|\nabla^k u\|_{X(\Sigma, \mu)}$ instead of $\|\, |\nabla^k u| \,\|_{X(\Sigma, \mu)}$.

We will also encounter Sobolev space $W^{1,1}(\Sigma, \mu)$, which is a weighted counterpart of the classical Sobolev space $W^{1,1}(\Sigma)$. We say that a function $u$ belongs to the space $W^{1,1}(\Sigma, \mu)$ if it is weakly differentiable in $\Sigma$, $u\in L^1(\Sigma, \mu)$ and $| \nabla u | \in L^1(\Sigma, \mu)$.

Now we turn back to the theory of re\-ar\-range\-ment\--in\-vari\-ant spaces. We present here their another important properties which we will need in what follows. 

To every re\-ar\-range\-ment\--in\-vari\-ant space $X$, there is associated another re\-ar\-range\-ment\--in\-vari\-ant space, which is related to its continuous dual, but which is usually more useful in the theory of Banach function spaces. The mapping $\rho'$ defined on $\mathcal{M}^{+}(R, \mu)$ by
\begin{align}\label{asno}
\rho'(g) = \sup\limits_{f \in \mathcal{M}^{+}(R, \mu), \rho(f) \leq 1}\int_{R} fg \, d\mu, \, \, \, \,  g \in \mathcal{M}^{+}(R, \mu),
\end{align}
is the \emph{associate norm} of the function norm $\rho$. We say that the space $X(\rho')$ is the associate space to the space $X(\rho)$ and we detone this space by $X'$.

For example, when $X = L^p(R, \mu)$ for $p\in[1, \infty]$, then $X' = L^{p'}(R, \mu)$. Here $p'\in[1, \infty]$ is the dual index defined by $\frac{1}{p} + \frac{1}{p'} = 1$.

An important property of Banach function spaces is that, if $X$ is a Banach function space, then $(X')' = X$.

For every Banach function norm $\rho$, \emph{the Hölder inequality}
\begin{align*}
\int_{R} \left\lvert fg \right\rvert \, d\mu \leq \rho(f) \rho'(g)
\end{align*}
holds for every $f, g \in \mathcal{M}(R, \mu)$.

For every $f \in \mathcal{M}(R, \mu)$ and $t \in (0, \infty)$, we have
\begin{align}\label{fco1}
\mu(\{x \in R; |f(x)| > f^{\ast}_{\mu}(t)\}) \leq t.
\end{align}
Furthermore, if $f^{\ast}_{\mu}(t) < \infty$ and $\mu(\{x \in R; |f(x)| > f^{\ast}_{\mu}(t) - \varepsilon\}) < \infty$ for some $\varepsilon > 0$, then
\begin{align}\label{fco2}
\mu(\{x \in R; |f(x)| \geq f^{\ast}_{\mu}(t)\}) \geq t. 
\end{align}

The \emph{fundamental function} $\varphi_{X} \colon [0, \mu(R)] \to [0, \infty]$ of the re\-ar\-range\-ment\--in\-vari\-ant space $X$ is the mapping defined by 
\begin{align*}
\varphi_{X}(t) = \left\lVert \chi_{E} \right\rVert_{X}, \,\,\,\, t \in [0, \mu(R)],
\end{align*}
where $E \subseteq R$ is an arbitrary set satisfying $\mu(E) = t$.
The definition is correct since if $E, F \subseteq R$ are sets such that $\mu(E) = \mu(F)$, then their characteristic functions $\chi_{E}$ and $\chi_{F}$ are equimeasurable.

Furthermore, if $g \in \mathcal{M}(R, \mu)$ and if $\rho$ is a re\-ar\-range\-ment\--in\-vari\-ant norm, we have
\begin{align}\label{assno}
\rho'(g) = \sup\limits_{f \in \mathcal{M}^{+}(R, \mu), \rho(f) \leq 1} \int_{R} f^{\ast}g^{\ast} \, d\mu.
\end{align}

\emph{The Hardy--Littlewood inequality} is very important in the theory of re\-ar\-range\-ment\--in\-vari\-ant spaces. It states that 
\begin{align*}
\int_{R} \left\lvert fg \right\rvert \, d\mu \leq \int\limits_{0}^{\infty} f^{\ast}(t) g^{\ast}(t) \, dt
\end{align*}
for every $f, g \in \mathcal{M}(R, \mu)$. In particular, by taking $g = \chi_E$, we have
\begin{align}\label{hli2}
\int_{E} \left\lvert f\right\rvert \, d\mu \leq \int\limits_{0}^{\mu(E)} f^{\ast}(t) \, dt
\end{align}
for each $\mu$-measurable set $E \subseteq R$.

Another important result in the theory of re\-ar\-range\-ment\--in\-vari\-ant spaces is the \emph{Hardy--Littlewood--P\'{o}lya principle}. For every re\-ar\-range\-ment\--in\-vari\-ant norm $\rho$, if $f, g \in \mathcal{M}^{+}(R, \mu)$ are such that
\begin{align*}
\int_{0}^{t} f^{\ast}(\tau) \, d\tau \leq \int_{0}^{t} g^{\ast}(\tau) \, d\tau
\end{align*}
for every $t \in (0, \infty)$, then
\begin{align}\label{klr}
\rho(f) \leq \rho(g).
\end{align}
We will also need the following fact. For every $t \in (0, \mu(R))$ and for every $f \in \mathcal{M}(R, \mu)$, we have
\begin{align}\label{ggt}
\int_{0}^{t}f^{\ast}(\tau) \; d\tau = \sup\left(\left\{ \int_{E} \left\lvert f\right\rvert \; d\mu; E \subseteq R, E \, \,  \text{$\mu$-measurable}, \mu(E) = t \right\}\right).
\end{align}

Each re\-ar\-range\-ment\--in\-vari\-ant space on $(R, \mu)$ can be represented as a re\-ar\-range\-ment\--in\-vari\-ant space on $(0, \mu(R))$. More precisely, if $X(R, \mu)$ is a re\-ar\-range\-ment\--in\-vari\-ant space, then there exists the unique re\-ar\-range\-ment\--in\-vari\-ant space $X(0, \mu(R))$ such that for every function $f \in X(R, \mu)$ it holds that
\begin{align*}
\left\lVert f \right\rVert_{X(R, \mu)} = \left\lVert f^{\ast}_{\mu} \right\rVert_{X(0, \mu(R))}.
\end{align*}
The re\-ar\-range\-ment\--in\-vari\-ant space $X(0, \mu(R))$ is called \emph{the representation space} of $X(R, \mu)$. For example, if $X(R, \mu) = L^p(R, \mu)$, then $X(0, \mu(R)) = L^p(0, \mu(R))$. 

On the other hand, for every $f \in \mathcal{M}(0, \mu(R))$, there exists a function $u \in \mathcal{M}(R, \mu)$ such that $f^{\ast}_{\lambda}(t) = u^{\ast}_{\mu}(t)$ for every $t \in (0, \infty)$.

Closely related to the nonincreasing rearrangement is \emph{the maximal nonincreasing rearrangement}. The \emph{maximal nonincreasing operator} 
\begin{align*}
P_{\mu} \colon \mathcal{M}(R, \mu) \rightarrow \mathcal{M}^{+}(0, \infty)
\end{align*}
is defined by
\begin{align*}
P_{\mu}(f)(t) = \frac{1}{t} \int_{0}^{t} f^{\ast}_{\mu}(\tau)\, d\tau, \, \, \, \, f \in \mathcal{M}(R, \mu), t \in (0, \infty).
\end{align*}
The image of a function $f \in \mathcal{M}(R, \mu)$ under the maximal nonincreasing operator $P_{\mu}$ is also commonly denoted by $f^{\ast \ast}_{\mu}$, and it is called \emph{the maximal nonincreasing function}. The maximal nonincreasing function is nonincreasing and we have $f^{\ast} \leq f^{\ast \ast}$.

If $X(0, \mu(R))$ is a re\-ar\-range\-ment\--in\-vari\-ant space and $h \in \mathcal{M}^{+}(0, \mu(R))$ is a nonincreasing function, we know thanks to (\ref{assno}) that
\begin{align*}
\left\lVert h \right\rVert_{X'(0, \mu(R))} = \sup\limits_{g \in \mathcal{M}^{+}(0, \mu(R)),  \left\lVert g \right\rVert_{X(0, \mu(R))} \leq 1}\int_{0}^{\mu(R)} h(t)g^{\ast}(t) \, dt.
\end{align*}
In general, when $h \in \mathcal{M}^{+}(0, \mu(R))$ is not necessarily nonincreasing, we only have
\begin{align}\label{pcps}
\left\lVert h \right\rVert_{X'(0, \mu(R))} \geq \sup\limits_{g \in \mathcal{M}^{+}(0, \mu(R)), \left\lVert g \right\rVert_{X(0, \mu(R))} \leq 1} \int_{0}^{\mu(R)} h(t)g^{\ast}(t) \, dt
\end{align}
owing to \eqref{asno}. However, it follows from \cite[Theorem 9.5]{CPS} and \cite[Theorem 3.10]{P1} that
\begin{align}\label{cpsp}
\left\lVert t^{\alpha}f^{\ast \ast}(t) \right\rVert_{X'(0, \mu(R))} \leq 4 \sup\limits_{g \in \mathcal{M}^{+}(0, \mu(R)), \left\lVert g \right\rVert_{X(0, \mu(R))} \leq 1}\int_{0}^{\mu(R)} t^{\alpha}f^{\ast \ast}(t)g^{\ast}(t) \, dt
\end{align}
for every $\alpha \in [0, 1]$ and $f\in \mathcal{M}(R, \mu)$. Inequalities (\ref{pcps}) and (\ref{cpsp}) mean that the norm of $t^{\alpha}f^{\ast \ast}(t)$ can be approached, up to a multiplicative constant, by nonincreasing functions in this case even though the function $t^{\alpha}f^{\ast \ast}(t)$ does not have to be nonincreasing.

We will continue by introducing the dilation operator. Let $\alpha \in (0, \infty)$. The \emph{dilation operator} $D_{\alpha} \colon \mathcal{M}^{+}(0, \infty) \rightarrow \mathcal{M}^{+}(0, \infty)$ is defined by $(D_{\alpha} f)(t) = f(\alpha t)$ for each $f \in \mathcal{M}^{+}(0, \infty)$ and $t \in (0, \infty)$. This operator is bounded on every re\-ar\-range\-ment\--in\-vari\-ant space over $(0, \infty)$. More precisely, there exists a constant $0 < C \leq \max\{1, \frac{1}{\alpha}\}$ such that 
\begin{align}\label{dila}
\left\lVert D_{\alpha} f \right\rVert_{X(0, \infty)} \leq C\left\lVert f \right\rVert_{X(0, \infty)}
\end{align}
for every $f \in \mathcal{M}^{+}(0, \infty)$, every $\alpha \in (0, \infty)$ and every re\-ar\-range\-ment\--in\-vari\-ant space $X(0, \infty)$.

We conclude this section by introduction the continuous embedding. We say that $X(R, \mu)$ is \emph{continuously embedded} into $Y(R, \mu)$ if for every function $u \in X(R,\mu)$ it holds that $u \in Y(R,\mu)$ and that $\left\lVert u \right\rVert_{Y(R, \mu)} \leq C \left\lVert u \right\rVert_{X(R, \mu)}$, where $C$ is a constant that does not depend on $u$. We denote the fact that $X(R, \mu)$ is continuously embedded into $Y(R, \mu)$ by $X(R, \mu) \hookrightarrow Y(R, \mu)$.
In fact, inclusion between Banach function spaces is always continuous in the sense that $X(R, \mu) \hookrightarrow Y(R, \mu)$ if and only if $X(R, \mu) \subseteq Y(R, \mu)$. If $X(R, \mu), Y(R, \mu)$ are Banach function spaces, it holds that $X(R, \mu) \hookrightarrow Y(R, \mu)$ if and only if $Y'(R, \mu) \hookrightarrow X'(R, \mu)$. 

In the rest of the paper, we will denote by $C, K, C_{i}, K_{i}$ positive finite constants whose exact values are not important for our purposes.

\section{Reduction principle and optimality}\label{sec:main}
\subsection{Reduction principle}\label{rp}
The goal of this subsection is to prove a suitable reduction principle. To prove it we need to derive a variant of the Pólya--Szeg\H{o} inequality. The prove of this theorem is at the end of this section.

\begin{prop}[Pólya--Szeg\H{o} inequality]\label{PS}
Let $X$ be a re\-ar\-range\-ment\--in\-vari\-ant space over $(\Sigma, \mu)$ and $u \in V^{1}_{0}X(\Sigma, \mu)$. Then $u^{\ast}_{\mu}$ is a locally absolutely continuous function on the interval $(0, \infty)$, and it holds that
\begin{align}\label{2.3}
\left\lVert t^{\frac{D - 1}{D}} \frac{d u^{\ast}_{\mu}}{dt}(t) \right\rVert_{X(0, \infty)} \leq C \left\lVert \nabla u \right\rVert_{X(\Sigma, \mu)},
\end{align}
where $C$ is a positive constant independent of $u$.
\end{prop}

In the remaining part of this section we prove the reduction principle. Recall that the parameteres $m$ and $D$ were introduced in \eqref{defD} and \eqref{defm}.

\begin{thm}[Reduction principle]\label{RP}
Let $X$ and $Y$ be re\-ar\-range\-ment\--in\-vari\-ant spaces over $(\Sigma, \mu)$. Then the following two statements are equivalent.
\begin{enumerate}
  \item For all functions $v \in V^{m}_{0}X(\Sigma, \mu)$ it holds that
\begin{align}\label{pnj3}
\left\lVert v \right\rVert_{Y(\Sigma, \mu)} \leq C_{1} \left\lVert \nabla^{m} v \right\rVert_{X(\Sigma, \mu)}.
\end{align}
  \item For all functions $f \in \mathcal{M}^{+}(0, \infty)$ it holds that
\begin{align}\label{scn}
\left\lVert \int_{t}^{\infty} f(\tau)\tau^{\frac{m}{D}-1} \; d\tau \right\rVert_{Y(0, \infty)} \leq C_{2} \left\lVert f \right\rVert_{X(0, \infty)}.
\end{align}
\end{enumerate}
Here, $C_{1}$ and $C_{2}$ are positive constants independent of $v$ and of $f$ respectively.
\end{thm}
The proof of this theorem will be divided into two steps, Proposition~\ref{d2} and Proposition~\ref{d3}.

\begin{rem}\label{trz3}
As an easy consequence of the definition of the associate norm, we obtain the fact that (\ref{scn}) is equivalent to the following assertion:
\begin{itemize}
\item[\emph{(2*)}]
\emph{For all functions $g \in \mathcal{M}^{+}(0, \infty)$ it holds that
\begin{align}\label{trd}
\left\lVert t^{\frac{m}{D}} g^{\ast \ast}(t) \right\rVert_{X'(0, \infty)} \leq C_{2} \left\lVert g \right\rVert_{Y'(0, \infty)}.
\end{align}}
\end{itemize}
\end{rem} 

\begin{rem}\label{gvlt}
We have
\begin{align*}
\begin{split}
& \left\lVert \chi_{(0, 1)}(t) \right\rVert_{Y(0, \infty)} \left\lVert t^{\frac{m}{D}-1} \chi_{(1, \infty)}(t) \right\rVert_{X'(0, \infty)} \\ & = \left\lVert \chi_{(0, 1)}(t) \right\rVert_{Y(0, \infty)} \sup\limits_{\left\lVert h\right\rVert_{X(0, \infty)}, h \geq 0} \int_{1}^{\infty}   \tau^{\frac{m}{D}-1} h(\tau) \,d\tau \\ & \leq \sup\limits_{\left\lVert h\right\rVert_{X(0, \infty)}, h \geq 0} \left\lVert \chi_{0, 1)}(t) \int_{t}^{\infty}   \tau^{\frac{m}{D}-1} h(\tau) \,d\tau \right\rVert_{Y(0, \infty)},
\end{split}
\end{align*}
so, \eqref{scn} implies that
\begin{align*}
t^{\frac{m}{D}-1} \chi_{(1, \infty)}(t) \in X'(0, \infty).
\end{align*}
\end{rem}

\begin{prop}\label{d2}
Let $X, Y$ be re\-ar\-range\-ment\--in\-vari\-ant spaces over $(\Sigma, \mu)$. Assume that there exists a positive constant $C_{1}$ such that for all functions $v \in V^{m}_{0}X(\Sigma, \mu)$ it holds that
\begin{align*}
\left\lVert v \right\rVert_{Y(\Sigma, \mu)} \leq C_{1} \left\lVert \nabla^{m} v \right\rVert_{X(\Sigma, \mu)}.
\end{align*}
Then there exists a positive constant $C_{2}$ such that for all functions $f \in \mathcal{M}^{+}(0, \infty)$ it holds that
\begin{align}\label{supp}
\left\lVert \int_{t}^{\infty} f(\tau)\tau^{\frac{m}{D}-1} \; d\tau \right\rVert_{Y(0, \infty)} \leq C_{2} \left\lVert f \right\rVert_{X(0, \infty)}. 
\end{align}
The constant $C_{2}$ depends only on the constant $C_{1}$, on $m$ and on $D$.
\end{prop}

Owing to this proposition the $n$-dimensional part of the reduction principle ((\ref{pnj3}) in Theorem~\ref{RP}) implies the one-dimensional part ((\ref{scn}) in Theorem~\ref{RP}). In the proof of the proposition, we need to use the following technical lemma, whose proof is straightforward and omitted.

\begin{lemma}\label{tcn}
Let $f \in \mathcal{M}^{+}(0, \infty) \cap L^{\infty}(0, \infty)$ be a function with a bounded support. Let $g \colon [0, \infty) \rightarrow [0, \infty)$ be the function defined by 
\begin{align}\label{gdef}
g(t) = \int_{t}^{\infty} f(\tau) \tau^{\frac{m}{D} - m} (\tau - t)^{m-1} \; d\tau, \; \; \; \; t \in [0, \infty).
\end{align}
Then $g \in \mathcal{C}^{m-1}(0, \infty)$ and
\begin{align*}
g^{(j)}(t) = (-1)^{j}\frac{(m-1)!}{(m-j-1)!} \int_{t}^{\infty} f(\tau) \tau^{\frac{m}{D} - m} (\tau - t)^{m-j-1} \; d\tau, \; \; \; \; t \in (0, \infty)
\end{align*}
for every $j \in \{1, \dots, m-1\}$. Moreover, $g^{(m-1)}$ is locally Lipschitz on $(0, \infty)$ and 
\begin{align*}
g^{(m)}(t) = (-1)^{m}(m-1)! f(t) t^{\frac{m}{D}-m}  
\end{align*}
for almost every $t \in (0, \infty)$.
\end{lemma}

We denote by $B_{\mu}$ the $\mu$-measure of the intersection of the unit ball in $\rn$ with $\Sigma$, i.e.,
\begin{align*}
B_{\mu} = \mu(\{x \in \Sigma; \left\lvert x \right\rvert \leq 1 \}).
\end{align*}

Now we prove Proposition~\ref{d2}.

\begin{myproof3}
Choose an arbitrary function $f \in \mathcal{M}^{+}(0, \infty)$. Observe that it is enough to prove the theorem for $\left\lVert f \right\rVert_{X(0, \infty)} < \infty$. Firstly, assume that the function $f$ belongs to $L^{\infty}(0, \infty)$ and that it has a bounded support. Define the function $v \colon \Sigma \rightarrow [0, \infty)$ by 
\begin{align}\label{fc1}
v(x) = \int_{\ballm \left\lvert x \right\rvert^{D}}^{\infty} f(\tau)\tau^{\frac{m}{D}-m}\left(\tau - \ballm \left\lvert x \right\rvert^{D}\right)^{m-1} \, d\tau, \, \, \, \, x \in \Sigma.
\end{align}
Now define the function $\sigma \colon \Sigma \rightarrow [0, \infty)$ by $\sigma(x) = \ballm \left\lvert x \right\rvert^{D}$, $x \in \Sigma$.
It holds that $\sigma \in \mathcal{C}^{\infty}(\Sigma)$. 
Clearly $v(x) = (g \circ \sigma)(x), x \in \Sigma$, where $g$ is the function from (\ref{gdef}), so $v \in \mathcal{C}^{m-1}(\Sigma)$ owing to Lemma~\ref{tcn} and the $m$th order weak partial derivatives exist (see \cite[Theorem 10.35]{L}) and are linear combinations of the functins \begin{align}\label{fc2}
x \mapsto \int_{\ballm \left\lvert x \right\rvert^{D}}^{\infty} f(\tau)\tau^{\frac{m}{D}-m}\left(\tau - \ballm \left\lvert x \right\rvert^{D}\right)^{m-l_{1}-1} \, d\tau \left\lvert x \right\rvert^{l_{1}(D-2) -2l_{2}} \prod\limits_{j=1}^{2(l_{1}+l_{2})-k} x_{i_{j}}
\end{align}
and
\begin{align}\label{f32}
x \mapsto f\left(\ballm \left\lvert x \right\rvert^{D}\right) \left\lvert x \right\rvert^{-m} \prod\limits_{j=1}^{m} x_{i_{j}} 
\end{align}
for the parametres $l_{1}, l_{2}$ satisfying $l_{1} \in \mathbb{N}, l_{1} \leq m-1, l_{2} \in \{0, \dots, m\}, 2(l_{1} + l_{2}) \geq m$. Since $f$ has a bounded support, we also have
\begin{align}\label{ww}
\left\lvert \nabla^{m} v(x) \right\rvert \leq K_{1}  \left(f\left(\ballm \left\lvert x \right\rvert^{D}\right) + \sum\limits_{l = 1}^{m-1} \int_{\ballm \left\lvert x \right\rvert^{D}}^{\infty} f(\tau) \tau^{\frac{m}{D}-l-1} \, d\tau \left\lvert x \right\rvert^{lD - m}\right)
\end{align}
for $\mu$-almost every $x \in \Sigma$ thanks to \eqref{fc2} and \eqref{f32}.

For every $l \in \{1, \dots m-1\}$, we now define the operator
\begin{align*}
F_{l} \colon (L^{1} + L^{\infty})(0, \infty) \rightarrow (L^{1} + L^{\infty})(0, \infty)
\end{align*}
by
\begin{align*}
F_{l}(\varphi)(t) =  t^{l-\frac{m}{D}}\int_{t}^{\infty} \varphi(\tau) \tau^{\frac{m}{D}-l-1} \, d\tau, \, \, t \in (0, \infty), \, \, \varphi \in (L^{1} + L^{\infty})(0, \infty).
\end{align*} 
For an arbitrary $l \in \{1, \dots m-1\}$ we have $\left\lVert F_{l} \right\rVert_{L^{\infty} \rightarrow L^{\infty}}  \leq \frac{D}{Dl-m}$ and $\left\lVert F_{l} \right\rVert_{L^{1} \rightarrow L^{1}} \leq \frac{D}{Dl-m+D}$. So, we obtain the fact that the operator $F_{l}$ is, owing to \cite[Chapter 3, Theorem 2.2]{BS}, bounded on the space $X(0, \infty)$.

Now define the functions $h \colon (0, \infty) \rightarrow [0, \infty)$ and $\omega \colon \Sigma \rightarrow [0, \infty)$. The function $h$ is defined by 
\begin{align*}
h(t) = f(t) + \sum_{l=1}^{m-1} F_{l}(f)(t), \, \, \, \, t \in (0, \infty).
\end{align*}
The function $\omega$ is defined by
\begin{align*}
\omega(x) = (h \circ \sigma)(x), \, \, \, \, x \in \Sigma. 
\end{align*}
Owing to (\ref{ww}) we obtain
\begin{align*}
\left\lvert \nabla^{m} v \right\rvert^{\ast}_{\mu}(t) \leq K_{2} \omega^{\ast}_{\mu}(t), \, \, \, \,   t \in (0, \infty).
\end{align*}
The functions $h$ and $\omega$ are equimeasurable. So, we have
\begin{align*}
\left\lvert \nabla^{m} v \right\rvert^{\ast}_{\mu}(t) \leq K_{2} h^{\ast}(t), \, \, \, \,   t \in (0, \infty),
\end{align*}
thanks to (\ref{rpp5}). Now we obtain
\begin{align}\label{bt}
\begin{split}
& \left\lVert \nabla^{m} v \right\rVert_{X(\Sigma, \mu)} \leq K_{2} \left\lVert h \right\rVert_{X(0, \infty)} \\ & \leq K_{2} \left(\left\lVert f \right\rVert_{X(0, \infty)} + \sum_{l = 1}^{m-1} \left\lVert F_{l}(f) \right\rVert_{X(0, \infty)}\right) \leq K_{3} \left\lVert f \right\rVert_{X(0, \infty)}.
\end{split}
\end{align}

From (\ref{bt}) it follows that $ v \in V^{m}X(\Sigma, \mu)$. Since the function $f$ has a bounded support, we obtain the fact that $v \in V^{m}_{0}X(\Sigma, \mu)$ owing to (\ref{fc1}). 

Now, we know that the functions $v$ and $g$ are equimeasurable since the mapping $\sigma$ is a measure-preserving transformation of the spaces $(\Sigma, \mu)$ and $([0, \infty), \lambda)$ (see \cite[Section 2.7]{BS}). We obtain the fact that
\begin{align}\label{fl1}
& \left\lVert v \right\rVert_{Y(\Sigma, \mu)} = \left\lVert g \right\rVert_{Y(0, \infty)}  \geq \left\lVert \int_{2t}^{\infty} f(\tau)\tau^{\frac{m}{D}-m}\left(\tau - t \right)^{m-1} \, d\tau \right\rVert_{Y(0, \infty)}  \\ & = 2^{1-m} \left\lVert \int_{2t}^{\infty} f(\tau)\tau^{\frac{m}{D}-m}\tau^{m-1} \, d\tau \right\rVert_{Y(0, \infty)} = 2^{1-m} \left\lVert \int_{2t}^{\infty} f(\tau)\tau^{\frac{m}{D}-1} \, d\tau \right\rVert_{Y(0, \infty)}. \notag
\end{align}

Finally, we have
\begin{align*}
\begin{split}
& \left\lVert \int_{t}^{\infty} f(\tau) \tau^{\frac{m}{D}-1} \, d\tau \right\rVert_{Y(0, \infty)} \leq 2 \left\lVert \int_{2t}^{\infty} f(\tau) \tau^{\frac{m}{D}-1}\, d\tau \right\rVert_{Y(0, \infty)} \leq 2^{m} \left\lVert v \right\rVert_{Y(\Sigma, \mu)} \\ & \leq 2^{m} C_{1} \left\lVert \nabla^{m} v \right\rVert_{X(\Sigma, \mu)} \leq 2^{m} C_{1} K_{3} \left\lVert f \right\rVert_{X(0, \infty)}.
\end{split}
\end{align*}
The first inequality holds by virtue of (\ref{dila}). The second inequality holds thanks to (\ref{fl1}). The last inequality holds owing to (\ref{bt}). So, we have proved the inequality (\ref{supp}) for bounded functions with bounded support. Now let $f$ be general. Define a sequence $\{f_{k}\}_{k=1}^{\infty}$ of~functions from~$\mathcal{M}^{+}(0, \infty)$ by $f_{k}(t) = \min\{f(t), k\}\chi_{(0, k)}(t), t \in (0, \infty), k \in \mathbb{N}$. Since (\ref{supp}) holds for every $f_{k}, k \in \mathbb{N}$, we obtain the fact that (\ref{supp}) also holds for $f$ thanks to the Fatou axiom of Banach function norms.
\end{myproof3}

Now we prove the remaining implication in Theorem~\ref{RP}.

\begin{prop}\label{d3}
Let $X, Y$ be re\-ar\-range\-ment\--in\-vari\-ant spaces over $(\Sigma, \mu)$. Assume that there exists a positive constant $C_{2}$ such that for all functions $f \in \mathcal{M}^{+}(0, \infty)$ it holds that
\begin{align}\label{rrv}
\left\lVert \int_{t}^{\infty} f(\tau)\tau^{\frac{m}{D}-1} \; d\tau \right\rVert_{Y(0, \infty)} \leq C_{2} \left\lVert f \right\rVert_{X(0, \infty)}.
\end{align}
Then there exists a positive constant $C_{1}$ such that for every $u \in V^{m}_{0}X(\Sigma, \mu)$ it holds that
\begin{align*}
\left\lVert u \right\rVert_{Y(\Sigma, \mu)} \leq C_{1} \left\lVert \nabla^{m} u \right\rVert_{X(\Sigma, \mu)}.
\end{align*}
The constant $C_{1}$ depends only on the constant $C_{2}$, on $m$ and on $D$.
\end{prop}
\begin{myproof}
The mapping $\sigma_{m, X} \colon \mathcal{M}^{+}(\Sigma, \mu) \rightarrow [0, \infty]$ defined by 
\begin{align*}
\sigma_{m, X}(v) = \left\lVert t^{\frac{m}{D}} v^{\ast \ast}_{\mu}(t) \right\rVert_{X'(0, \infty)}, \, \, \, \, v \in \mathcal{M}^{+}(\Sigma, \mu),
\end{align*}
is a re\-ar\-range\-ment\--in\-vari\-ant Banach function norm if and only if 
\begin{align}\label{ghg}
t^{\frac{m}{D}-1} \chi_{(1, \infty)}(t) \in X'(0, \infty)
\end{align}
(see \cite[Theorem 5.4]{CPS} and \cite[Theorem 4.4]{EMMP}).
Since (\ref{rrv}) holds we have the fact that \eqref{ghg} is true thanks to Remark~\ref{gvlt}. Hence $\sigma_{m, X}$ is a re\-ar\-range\-ment\--in\-vari\-ant Banach function norm. We denote the respective space by $Z_{m}X(\Sigma, \mu)$.

We want to prove that for every function $u \in V^{m}_{0}X(\Sigma, \mu)$ it holds that
\begin{align}\label{rlst}
\left\lVert u \right\rVert_{Y(\Sigma, \mu)} \leq C_{2} \left\lVert u \right\rVert_{Z'_{m}X(\Sigma, \mu)} \leq C_{2}K_{m} \left\lVert \nabla^{m} u \right\rVert_{X(\Sigma, \mu)},
\end{align}
where $K_{m}$ is a positive constant. The first inequality follows from (\ref{trd}). We prove the second inequality by induction on $m$. Firstly, we assume that $m = 1$. Choose an arbitrary function $u \in V_{0}^{1}X(\Sigma, \mu)$. We have
\begin{align*}
\begin{split}
& \left\lVert u \right\rVert_{Z'_{1}X(\Sigma, \mu)} = \left\lVert u^{\ast}_{\mu} \right\rVert_{Z'_{1}X(0, \infty)} = \left\lVert - \int_{t}^{\infty} \frac{d u^{\ast}_{\mu}}{d\tau}(\tau) \, d\tau \right\rVert_{Z'_{1}X(0, \infty)} \\ & = \left\lVert \int_{t}^{\infty} \left ( \tau^{\frac{D-1}{D}} \frac{d u^{\ast}_{\mu}}{d\tau}(\tau) \right ) \tau^{\frac{1-D}{D}} \, d\tau  \right\rVert_{Z'_{1}X(0, \infty)} \leq \left\lVert t^{\frac{D-1}{D}} \frac{d u^{\ast}_{\mu}}{d t}(t) \right\rVert_{X(0, \infty)} \\ & \leq C_{3} \left\lVert \nabla u \right\rVert_{X(\Sigma, \mu)}.
\end{split}
\end{align*}
The second equality is true owing to the fact that $u_{\mu}$ is locally absolutely continuous on $(0, \infty)$ (see Proposition~\ref{PS}). The first inequality holds by virtue of (\ref{trd}) and the second one is true thanks to the Pólya--Szeg\H{o} inequality (Proposition~\ref{PS}). So, we proved \eqref{rlst} for $m = 1$.

Now assume that $1 < m < D$ is arbitrary and that \eqref{rlst} holds for $m-1$. Choose an arbitrary function $u \in V_{0}^{m}X(\Sigma, \mu)$ and $i \in \{1, \dots, n\}$. Then the weak partial derivative $\frac{\partial u}{\partial x_{i}}$ belongs to $V_{0}^{m-1}X(\Sigma, \mu)$. So, we can use the induction hypothesis to obtain 
\begin{align*}
\left\lVert \frac{\partial u}{\partial x_{i}} \right\rVert_{Z'_{m-1}X(\Sigma, \mu)} \leq K_{m-1} \left\lVert  \nabla^{m-1} \frac{\partial u}{\partial x_{i}} \right\rVert_{X(\Sigma, \mu)} \leq K_{m-1} \left\lVert \nabla^{m} u \right\rVert_{X(\Sigma, \mu)}.
\end{align*}
It means that 
\begin{align}\label{ddd}
\left\lVert \nabla u \right\rVert_{Z'_{m-1}X(\Sigma, \mu)} \leq n K_{m-1} \left\lVert \nabla^{m} u \right\rVert_{X(\Sigma, \mu)} < \infty.
\end{align}
It follows that $u$ belongs to $V_{0}^{1}Z'_{m-1}X(\Sigma, \mu)$. Now it can be straightforwardly proved that
\begin{align*}
t^{\frac{1}{D}-1}\chi_{(1, \infty)}(t) \in Z_{m-1}X(0, \infty),
\end{align*}
(for details see \cite[Theorem 2.2]{M}), which means that the mapping $u \mapsto \left\lVert t^{\frac{1}{D}} u^{\ast \ast}_{\mu_{D}}(t) \right\rVert_{Z_{m-1}X(0, \infty)} = \left\lVert u \right\rVert_{Z_{1}(Z_{m-1}X)(0, \infty)}$ is a re\-ar\-range\-ment\--in\-vari\-ant Banach function norm (cf. (\ref{ghg})). From the case $m = 1$ it now follows that
\begin{align}\label{e2e}
\left\lVert u \right\rVert_{(Z_{1}(Z_{m-1}X))'(\Sigma, \mu)} \leq C_{3} \left\lVert \nabla u \right\rVert_{Z'_{m-1}X(\Sigma, \mu)}.
\end{align}
Owing to \cite[Theorem~3.4]{CCP} and \cite[Proposition~5.1]{M2} (cf. \cite[Theorem 9.5]{CPS}), we obtain
\begin{align*}
\begin{split}
& \left\lVert v \right\rVert_{Z_{1}(Z_{m-1}X)(\Sigma, \mu)} = \left\lVert t^{\frac{m-1}{D}} \left(\tau^{\frac{1}{D}} v^{\ast \ast}_{\mu}(\tau)\right)^{\ast \ast}(t) \right\rVert_{X'(0, \infty)} \\ & \leq C_{4} \left\lVert t^{\frac{m}{D}} v^{\ast \ast}_{\mu}(t) \right\rVert_{X'(0, \infty)} = C_{4} \left\lVert v \right\rVert_{Z_{m}X(\Sigma, \mu)}
\end{split}
\end{align*}
for every $v \in \mathcal{M}^{+}(\Sigma, \mu)$. From the previous inequality we obtain
\begin{align}\label{ffff}
\left\lVert v \right\rVert_{Z'_{m}X(\Sigma, \mu)} \leq C_ {5} \left\lVert v \right\rVert_{(Z_{1}(Z_{m-1}X))'(\Sigma, \mu)}
\end{align}
for every $v \in \mathcal{M}^{+}(\Sigma, \mu)$. By virtue of (\ref{ddd}), (\ref{e2e}) and (\ref{ffff}), we obtain the fact that (\ref{rlst}) holds for $m$.
\end{myproof}

\subsection{Optimality}\label{opt}
In this subsection we will describe the smallest target space in the inequality (\ref{pnj3}) among all spaces $Y(\Sigma, \mu)$ for a given space $X(\Sigma, \mu)$ and the largest domain space  in the inequality (\ref{pnj3}) among all spaces $X(\Sigma, \mu)$ for a given space $Y(\Sigma, \mu)$. 

Firstly we describe the optimal target space in the following theorem.
\begin{thm}\label{zopt}
Let $X$ be a re\-ar\-range\-ment\--in\-vari\-ant space over $(\Sigma, \mu)$. Assume that 
\begin{align}\label{jnvm}
t^{\frac{m}{D}-1}\chi_{(1, \infty)}(t) \in X'(0, \infty).
\end{align}
Then there exists a positive constant $C$, which depends only on $m$ and on $D$, such that
\begin{align}\label{iqsz}
\left\lVert u \right\rVert_{Z'_{m}X(\Sigma, \mu)} \leq C \left\lVert \nabla^{m} u \right\rVert_{X(\Sigma, \mu)}, \, \, \, \, u \in V^{m}_{0}X(\Sigma, \mu).
\end{align}
Moreover, the space $Z'_{m}X(\Sigma, \mu)$ is the optimal space in the previous inequality among all re\-ar\-range\-ment\--in\-vari\-ant spaces in the following way. If $Y(\Sigma, \mu)$ is a re\-ar\-range\-ment\--in\-vari\-ant space satisfying 
\begin{align}\label{iqsy}
\left\lVert u \right\rVert_{Y(\Sigma, \mu)} \leq \widetilde{C} \left\lVert \nabla^{m} u \right\rVert_{X(\Sigma, \mu)}, \, \, \, \, u \in V^{m}_{0}X(\Sigma, \mu),
\end{align}
with a positive constant $\widetilde{C}$ that does not depend on $u$, then 
\begin{align}\label{ebb}
Z'_{m}X(\Sigma, \mu) \hookrightarrow Y(\Sigma, \mu).
\end{align}
On the other hand, if (\ref{jnvm}) is not true, then the inequality (\ref{iqsy}) does not hold for any re\-ar\-range\-ment\--in\-vari\-ant space $Y(\Sigma, \mu)$.
\end{thm}
\begin{myproof}
The inequality (\ref{iqsz}) follows from the proof of Proposition~\ref{d3}. Assume that (\ref{jnvm}) holds and that $Y(\Sigma, \mu)$ is a re\-ar\-range\-ment\--in\-vari\-ant space satisfying (\ref{iqsy}). Owing to Proposition~\ref{d2} we obtain the fact that $Y(\Sigma, \mu)$ satisfies the inequality (\ref{scn}). It means that we can use the first inequality in (\ref{rlst}) to obtain the embedding (\ref{ebb}). 

Now, assume that (\ref{jnvm}) is not true. Then there does not exist any re\-ar\-range\-ment\--in\-vari\-ant space $Y(\Sigma, \mu)$ such that (\ref{scn}) holds thanks to Remark~\ref{gvlt}. Then by Proposition~\ref{d2} there is not any re\-ar\-range\-ment\--in\-vari\-ant space $Y(\Sigma, \mu)$ such that (\ref{iqsy}) is true.
\end{myproof}

Let $Y$ be a re\-ar\-range\-ment\--in\-vari\-ant space over $(\Sigma, \mu)$. Assume that 
\begin{align}\label{jnvm22}
\inf\limits_{t \in [1, \infty)} \frac{t^{1-\frac{m}{D}}}{\varphi_{Y}(t)} > 0,
\end{align}
where $\varphi_{Y}$ is the fundamental function of the space $Y(\Sigma, \mu)$.Then the mapping $\nu_{m, Y} \colon \mathcal{M}^{+}(\Sigma, \mu) \rightarrow [0, \infty]$ defined by
\begin{align*}
\nu_{m, Y}(u) = \sup\limits_{v \sim u, v \geq 0} \left\lVert \int_{t}^{\infty} v(\tau) \tau^{\frac{m}{D} -1} \, d\tau \right\rVert_{Y(0, \infty)},
\end{align*}
where the supremum is taken over all functions $v \in \mathcal{M}^{+}(0, \infty)$ equimeasurable with $u$, is a re\-ar\-range\-ment\--in\-vari\-ant Banach function norm (for the proof see \cite[Theorem 4.1]{EMMP}). The respective space we will denote by $U_{m}Y(\Sigma, \mu)$. 

Now we can finally describe the largest domain space.

\begin{thm}\label{zopt22}
Let $Y$ be a re\-ar\-range\-ment\--in\-vari\-ant space over $(\Sigma, \mu)$. Assume that the condition (\ref{jnvm22}) is satisfied. Then there exists a positive constant $C$, which depends only on $m$ and on $D$, such that
\begin{align}\label{iqsz22}
\left\lVert u \right\rVert_{Y(\Sigma, \mu)} \leq C \left\lVert \nabla^{m} u \right\rVert_{U_{m}Y(\Sigma, \mu)}, \, \, \, \, u \in V^{m}_{0}U_{m}Y(\Sigma, \mu).
\end{align}
Moreover, the space $U_{m}Y(\Sigma, \mu)$ is the optimal space in the previous inequality among all re\-ar\-range\-ment\--in\-vari\-ant spaces in the following way. If $X(\Sigma, \mu)$ is a re\-ar\-range\-ment\--in\-vari\-ant space satisfying 
\begin{align}\label{iqsy22}
\left\lVert u \right\rVert_{Y(\Sigma, \mu)} \leq \widetilde{C} \left\lVert \nabla^{m} u \right\rVert_{X(\Sigma, \mu)}, \, \, \, \, u \in V^{m}_{0}X(\Sigma, \mu),
\end{align}
with a positive constant $\widetilde{C}$ that does not depend on $u$, then 
\begin{align}\label{ebb22}
X(\Sigma, \mu) \hookrightarrow U_{m}Y(\Sigma, \mu).
\end{align}
On the other hand, if (\ref{jnvm22}) is not true, then the inequality (\ref{iqsy22}) does not hold for any re\-ar\-range\-ment\--in\-vari\-ant space $X(\Sigma, \mu)$.
\end{thm}

\begin{myproof}
Let $f \in \mathcal{M}^{+}(0, \infty)$. Then 
\begin{align*}
\begin{split}
& \left\lVert \int_{t}^{\infty} f(\tau) \tau^{\frac{m}{D} -1} \, d\tau \right\rVert_{Y(0, \infty)} \leq \sup\limits_{g \sim f, g \geq 0} \left\lVert \int_{t}^{\infty} g(\tau) \tau^{\frac{m}{D} -1} \, d\tau \right\rVert_{Y(0, \infty)} \\ & = \left\lVert f \right\rVert_{U_{m}Y(0, \infty)},
\end{split}
\end{align*}
so, the inequality (\ref{iqsz22}) follows from Proposition~\ref{d3}. 

Now, assume that $X(\Sigma, \mu)$ is a re\-ar\-range\-ment\--in\-vari\-ant space satisfying (\ref{iqsy22}) and let  $f, g \in \mathcal{M}^{+}(0, \infty)$ be equimeasurable. Then, by virtue of Proposition~\ref{d2}, we obtain  
\begin{align*}
\begin{split}
& \left\lVert f \right\rVert_{U_{m}Y(0, \infty)} = \sup\limits_{g \sim f, g \geq 0} \left\lVert \int_{t}^{\infty} g(\tau) \tau^{\frac{m}{D} -1} \, d\tau \right\rVert_{Y(0, \infty)} \leq K \sup\limits_{g \sim f, g \geq 0} \left\lVert g \right\rVert_{X(0, \infty)} \\ & = K \left\lVert f \right\rVert_{X(0, \infty)}.
\end{split}
\end{align*}
It means that (\ref{ebb22}) is true.

On the other hand, if the condition (\ref{jnvm22}) is not satisfied, then it can be proved analogously as in the proof of \cite[Theorem 4.1]{EMMP} that for each $u \in \mathcal{M}^{+}(\Sigma, \mu)$ such that $u^{\ast}_{\mu} = \chi_{(0, 1)}$ it holds that
\begin{align*}
\nu_{m, Y}(u) = \sup\limits_{v \sim u, v \geq 0} \left\lVert \int_{t}^{\infty} v(\tau) \tau^{\frac{m}{D} -1} \, d\tau \right\rVert_{Y(0, \infty)} = \infty.
\end{align*}
By the nontriviality axiom of the Banach function norm, we know that $u$ belongs to an arbitrary re\-ar\-range\-ment\--in\-vari\-ant space $X(\Sigma, \mu)$. So, we obtain the fact that \eqref{scn} does not hold for any re\-ar\-range\-ment\--in\-vari\-ant space $X(\Sigma, \mu)$. It means that the inequality (\ref{iqsy22}) does not hold for any re\-ar\-range\-ment\--in\-vari\-ant space $X(\Sigma, \mu)$ owing to the reduction principle (Theorem~\ref{RP}).
\end{myproof}

At the end of the section we prove Proposition \ref{PS}. The proof is based on the proofs of theorems \cite[Lemma 4.1]{CP1}, \cite[Lemma 3.3]{CP2} and \cite[Lemma 1.E]{T}.

\begin{myproof4}
Firstly, we prove the proposition for nonnegative $u$. We start with the proof of the local absolute continuity of the function $u^{\ast}_{\mu}$. Let $\{(a_{m}, b_{m})\}_{m\in M}$ be a countable system of pairwise disjoint nonempty bounded intervals. For each $m \in M$ define the function $f_{m} \colon \mathbb{R} \rightarrow \mathbb{R}$ in the following way: \[ f_{m}(y) = \left \{
\begin{array}{ll} 0 & \text{if} \, \,  y \leq u^{\ast}_{\mu}(b_{m}), \\
                       y - u^{\ast}_{\mu}(b_{m}) & \text{if} \, \, u^{\ast}_{\mu}(b_{m}) < y < u^{\ast}_{\mu}(a_{m}), \\
                       u^{\ast}_{\mu}(a_{m}) - u^{\ast}_{\mu}(b_{m}) & \text{if} \, \, u^{\ast}_{\mu}(a_{m}) \leq y. \\
\end{array}
\right. \]
For each $m \in M$ we now set $v_{m} = f_{m} \circ u$. It also holds that
\begin{align}\label{rr} 
\mu(\{x \in \Sigma; u(x) > u^{\ast}_{\mu}(b_{m})\}) < \infty
\end{align}
since $u \in V^{1}_{0}X(\Sigma, \mu)$. So, the function $v_{m}$ is bounded and can be nonzero in a set of finite $\mu$-measure only. We obtain the fact that $v_{m} \in L^{1}(\Sigma, \mu)$. Now we use the chain rule for Sobolev functions to obtain the fact that $v_{m}$ is weakly differentiable in $\Sigma$ and $\nabla v_{m} = \nabla u \chi_{\{u^{\ast}_{\mu}(b_{m}) < u < u^{\ast}_{\mu}(a_{m})\}}$ $\mu$-almost everywhere in $\Sigma$. From this equality we get $\left\lvert \nabla v_{m} \right\rvert = 0$ $\mu$-almost everywhere in the set $\{x \in\Sigma; u(x) \leq u^{\ast}_{\mu}(b_{m})\}$. We know that the function $\nabla u \in X(\Sigma, \mu)$. So, by virtue of (\ref{emb}) it follows that $\nabla u \in L^{1}(E, \mu)$ for every $\mu$-measurable set $E \subseteq \Sigma, \mu(E) < \infty$. We can now again use (\ref{rr}) to obtain $\nabla v_{m} \in L^{1}(\Sigma, \mu)$.

Now we can use the coarea formula (see \cite{MSZ}, \cite{S:07}) for the functions $v_{m}$, $m \in M$, and the isoperimetric inequality (see \cite[Theorem 1.3]{CR-OS}). We obtain 
\begin{align}\label{CI}
\begin{split}
& \int_{\bigcup_{m \in M} \{u^{\ast}_{\mu}(b_{m}) < u < u^{\ast}_{\mu}(a_{m})\}} \left\lvert \nabla u(x) \right\rvert \, d\mu(x) \\ & \geq C_{iso} \sum\limits_{m \in M} \int^{u^{\ast}_{\mu}(a_{m}) - u^{\ast}_{\mu}(b_{m})}_{0} \mu(\{x \in \Sigma ; v_{m}(x) > t\})^{\frac{D-1}{D}} \, dt \\ & = C_{iso} \sum\limits_{m \in M} \int^{u^{\ast}_{\mu}(a_{m})}_{u^{\ast}_{\mu}(b_{m})}  \mu(\{ x \in \Sigma; u(x) > t \})^{\frac{D-1}{D}} \, dt,
\end{split}
\end{align}
where $C_{iso}$ is the isoperimetric constant. Now we derive an upper estimate of
\begin{align*}
\int_{\bigcup_{m \in M} \{u^{\ast}_{\mu}(b_{m}) < u < u^{\ast}_{\mu}(a_{m})\}} \left\lvert \nabla u(x) \right\rvert \, d\mu(x).
\end{align*}
We obtain the fact that
\begin{align}\label{x}
\begin{split}
& \int_{\bigcup_{m \in M} \{u^{\ast}_{\mu}(b_{m}) < u < u^{\ast}_{\mu}(a_{m})\}} \left\lvert \nabla u(x) \right\rvert \, d\mu(x) \\ & \leq \int^{\sum\limits_{m \in M} \mu(\{u^{\ast}_{\mu}(b_{m}) < u < u^{\ast}_{\mu}(a_{m})\})}_{0} \left\lvert \nabla u \right\rvert^{\ast}_{\mu}(t) \, dt \leq \int^{\sum\limits_{m \in M} (b_{m} - a_{m})}_{0} \left\lvert \nabla u \right\rvert^{\ast}_{\mu}(t) \, dt.
\end{split}
\end{align}
The first inequality holds by virtue of the Hardy–-Littlewood inequality (\ref{hli2}). We can verify the last inequality in the following way. We have 
\begin{align*}
\begin{split}
& \mu(\{u_{\mu}^{\ast}(b_{m}) < u < \mu_{\mu}^{\ast}(a_{m})\}) = \mu(\{u_{\mu}^{\ast}(b_{m}) < u \}) - \mu(\{u_{\mu}^{\ast}(a_{m}) \leq u\}) \\ & \leq b_{m} - a_{m},
\end{split}
\end{align*}
where we used (\ref{fco1}) and (\ref{fco2}).

In the following part of the proof, we will assume that all the intervals $(a_{m}, b_{m}), m \in M$, are contained in an interval $[a, b] \subseteq (0, \infty)$. We prove that the function $u^{\ast}_{\mu}$ is absolutely continuous on the interval $[a, b]$. We can assume that $u_{\mu}^{\ast}(a) > 0$. We set
\begin{align*}
K = \mu(\{ x \in \Sigma; u(x) \geq u^{\ast}_{\mu}(a) \}).
\end{align*}
Then we have $K < \infty$ since $u \in V^{1}_{0}X(\Sigma, \mu)$. Since $u_{\mu}^{\ast}(a) > 0$, we can use (\ref{fco2}) to obtain $K > 0$. Owing to (\ref{CI}) we obtain the fact that
\begin{align}\label{M}
\begin{split}
& \int_{\bigcup_{m \in M} \{u^{\ast}_{\mu}(b_{m}) < u < u^{\ast}_{\mu}(a_{m})\}} \left\lvert \nabla u(x) \right\rvert \, d\mu(x) \\ & \geq C_{iso}\sum\limits_{m \in M} \int_{u^{\ast}_{\mu}(b_{m})}^{u^{\ast}_{\mu}(a_{m})} \mu(\{x \in \Sigma; u(x) \geq u^{\ast}_{\mu}(a)\})^{\frac{D-1}{D}} \, dt \\
& = C_{iso} K^{\frac{D-1}{D}} \sum\limits_{m \in M} (u^{\ast}_{\mu}(a_{m}) - u^{\ast}_{\mu}(b_{m})).
\end{split}
\end{align}
It follows that
\begin{align}\label{F}
\begin{split}
& \sum\limits_{m \in M} (u^{\ast}_{\mu}(a_{m}) - u^{\ast}_{\mu}(b_{m})) \\ & \leq C_{iso}^{-1} K^{\frac{1-D}{D}} \int_{\bigcup_{m \in M} \{u^{\ast}_{\mu}(b_{m}) < u < u^{\ast}_{\mu}(a_{m})\}} \left\lvert \nabla u(x) \right\rvert \, d\mu(x) \\ & \leq  C_{iso}^{-1} K^{\frac{1-D}{D}} \int^{\sum\limits_{m \in M} (b_{m} - a_{m})}_{0} \left\lvert \nabla u \right\rvert^{\ast}_{\mu}(t) \, dt.
\end{split}
\end{align}
The first inequality holds due to (\ref{M}). The second inequality holds by virtue of (\ref{x}). 

Next we want to prove that
\begin{align}\label{int}
 \int_{0}^{t} \left\lvert \nabla u \right\rvert^{\ast}_{\mu} (\tau) \, d\tau < \infty
\end{align}
for every $t \in (0, \infty)$. Since $\left\lvert \nabla u \right\rvert \in (L^{1} + L^{\infty})(\Sigma, \mu)$, it follows that $\left\lvert \nabla u \right\rvert^{\ast}_{\mu} \in L^{1}(0, t)$. So, (\ref{int}) is true.

Since we know that $\left\lvert \nabla u \right\rvert^{\ast}_{\mu}$ is integrable over an arbitrary bounded interval $(0, t)$, we can use (\ref{F}) to obtain the fact that the function $u^{\ast}_{\mu}$ is absolutely continuous on the interval $[a, b]$. It follows that it is locally absolutely continuous on the interval $(0, \infty)$, which is the desired result. 

It remains to prove the inequality (\ref{2.3}). From now we do not anymore assume that the intervals $(a_{m}, b_{m})$ are contained in $[a, b]$. Define the function $\phi \colon (0, \infty) \rightarrow [0, \infty)$ by $\phi(t) = - C_{iso} t^{\frac{D-1}{D}} \frac{du^{\ast}_{\mu}}{dt}(t), \, \, t \in (0, \infty)$. We show that
\begin{equation}\label{star}
\int_{0}^{t}\phi^{\ast}(\tau) \; d\tau \leq \int_{0}^{t}\left\lvert \nabla u \right\rvert^{\ast}_{\mu} (\tau) \; d\tau,  \, \, \, \, t \in (0, \infty).
\end{equation}
Choose $t \in (0, \infty)$ arbitrarily. By virtue of (\ref{ggt}), we know that it is enough to prove that for every measurable set $E \subseteq (0, \infty)$ such that $\lambda(E) = t$, it holds that
\begin{align}\label{act}
\int_{E} \phi(\tau)\; d\tau \leq \int_{0}^{t}\left\lvert \nabla u \right\rvert^{\ast}_{\mu} (\tau) \; d\tau.
\end{align}
Now choose an arbitrary $m \in M$. We have
\begin{align}\label{b}
\int_{a_{m}}^{b_{m}} \phi(\tau) \, d\tau = - \int_{a_{m}}^{b_{m}} C_{iso} \tau^{\frac{D-1}{D}}\frac{du^{\ast}_{\mu}}{d\tau}(\tau) \, d\tau.
\end{align}
Since the function $u^{\ast}_{\mu}$ is absolutely continuous and nonincreasing on the interval $[a_{m}, b_{m}]$, and the function $\tau \mapsto \mu(\{x \in \Sigma; u(x) > \tau\})^{\frac{D-1}{D}}, \tau \in (0, \infty),$ is nonnegative on the interval $(a_{m}, b_{m})$, we obtain the fact that (see \cite[page 156]{R})
\begin{align}\label{c}
\begin{split}
& \int_{u^{\ast}_{\mu}(b_{m})}^{u^{\ast}_{\mu}(a_{m})} C_{iso} \mu(\{x \in \Sigma; u(x) > s \})^{\frac{D-1}{D}} \; ds \\ & = -\int_{a_{m}}^{b_{m}} C_{iso} \mu(\{x \in \Sigma; u(x) > u^{\ast}_{\mu}(\tau)\})^{\frac{D-1}{D}} \frac{du^{\ast}_{\mu}}{d\tau}(\tau) \; d\tau.
\end{split}
\end{align}
Now we observe that
\begin{align}\label{eq}
\begin{split}
& \int_{a_{m}}^{b_{m}} C_{iso} \mu(\{x \in \Sigma; u(x) > u^{\ast}_{\mu}(\tau)\})^{\frac{D-1}{D}} \frac{du^{\ast}_{\mu}}{d\tau}(\tau) \; d\tau \\ & = \int_{a_{m}}^{b_{m}} C_{iso} \tau^{\frac{D-1}{D}} \frac{du^{\ast}_{\mu}}{d\tau}(\tau) \; d\tau.
\end{split}
\end{align}
Owing to (\ref{fco1}) we know that for every $\tau \in (a_{m}, b_{m})$ it holds that 
\begin{align*}
\mu(\{x \in \Sigma; u(x) > u^{\ast}_{\mu}(\tau)\}) \leq \tau.
\end{align*}
Choose an arbitrary $\tau \in (a_{m}, b_{m})$ such that the function $u^{\ast}_{\mu}$ is differentiable at $\tau$. Assume that 
\begin{align*}
\mu(\{x \in \Sigma; u(x) > u^{\ast}_{\mu}(\tau)\}) < \tau.
\end{align*}
Then there exists $\delta > 0$, such that for every $s \in (\tau - \delta, \tau)$ it holds that
\begin{align*}
\mu(\{x \in \Sigma; u(x) > u^{\ast}_{\mu}(\tau)\}) < s.
\end{align*}
From this inequality we obtain the fact that for each $s \in (\tau - \delta, \tau)$
\begin{align*}
\begin{split}
& u^{\ast}_{\mu}(s) = \inf\left\{\alpha > 0; \mu(\{x \in \Sigma; u(x) > \alpha\}) \leq s\right\} \\ & \leq u^{\ast}_{\mu}(\tau).
\end{split}
\end{align*}
Because the function~$u^{\ast}_{\mu}$ is nonincreasing it holds that $u^{\ast}_{\mu}(s) = u^{\ast}_{\mu}(\tau)$ for all $s \in (\tau - \delta, \tau)$. Since the function $u^{\ast}_{\mu}$ is differentiable at $\tau$, we have $\frac{du^{\ast}_{\mu}}{d\tau}(\tau) = 0$. So, we have just proved (\ref{eq}).
By (\ref{b}), (\ref{c}) and (\ref{eq}) we have
\begin{align}\label{e}
\int_{a_{m}}^{b_{m}} \phi(\tau) \, d\tau = \int_{u^{\ast}_{\mu}(b_{m})}^{u^{\ast}_{\mu}(a_{m})} C_{iso} \mu(\{x \in \Sigma; u(x) > \tau \})^{\frac{D-1}{D}} \; d\tau.
\end{align}
Now we obtain
\begin{align}\label{fin}
\begin{split}
& \int_{\bigcup_{m \in M}(a_{m}, b_{m})} \phi(\tau) \, d\tau = \sum\limits_{m \in M} \int_{u^{\ast}_{\mu}(b_{m})}^{u^{\ast}_{\mu}(a_{m})} C_{iso} \mu(\{x \in \Sigma; u(x) > \tau \})^{\frac{D-1}{D}} \; d\tau \\ & \leq \int_{\bigcup_{m \in M} \{u^{\ast}_{\mu}(b_{m}) < u < u^{\ast}_{\mu}(a_{m})\}} \left\lvert \nabla u(x) \right\rvert \, d\mu(x) \\ & \leq \int^{\sum\limits_{m \in M} (b_{m} - a_{m})}_{0} \left\lvert \nabla u \right\rvert^{\ast}_{\mu}(\tau) \, d\tau.
\end{split}
\end{align}
The equality holds thanks to (\ref{e}) and to the fact that all the intervals $(a_{m}, b_{m}), m \in M$, are pairwise disjoint. The first inequality holds thanks to (\ref{CI}). The second inequality is true by virtue of (\ref{x}). Now choose an arbitrary measurable set $E \subseteq (0, \infty), \lambda(E) = t$. Then for every $\varepsilon > 0$, there exists a countable system $\{(a_{m}, b_{m})\}_{m \in M}$ of pairwise disjoint nonempty bounded intervals such that
\begin{align*}
E \subseteq \bigcup_{m \in M} (a_{m}, b_{m}) \; \; \; \text{and} \; \; \; \lambda\left(\bigcup_{m \in M} (a_{m}, b_{m}) \setminus E \right) < \varepsilon.
\end{align*}
So, choose an arbitrary $\varepsilon > 0$, and let $\{(a_{m}, b_{m})\}_{m \in M}$ be such a system of intervals. From (\ref{fin}) we obtain
\begin{align*}
\int_{E} \phi(\tau)\, d\tau \leq \int_{\bigcup_{m \in M}(a_{m}, b_{m})} \phi(\tau) \, d\tau \leq \int^{t + \varepsilon}_{0} \left\lvert \nabla u \right\rvert^{\ast}_{\mu}(\tau) \, d\tau.
\end{align*}
We already know from (\ref{int}) that for all $s \in (0, \infty)$ the function $\left\lvert \nabla u \right\rvert^{\ast}_{\mu}$ is integrable over $(0, s)$. So, we obtain
\begin{align*}
\int_{E} \phi(\tau) \, d\tau \leq \int^{t}_{0} \left\lvert \nabla u \right\rvert^{\ast}_{\mu}(\tau) \, d\tau.
\end{align*}
This means that (\ref{act}) is true. So, we have just proved (\ref{star}). The inequality (\ref{2.3}) now follows from the Hardy--Littlewood--Pólya principle (\ref{klr}). 

We have proved the proposition for nonnegative $u$. Now let $u \in V_{0}^{1}X(\Sigma, \mu)$ be general. It is true that $\left\lvert u \right\rvert \in V_{0}^{1}X(\Sigma, \mu)$. So, the proposition holds for $\left\lvert u \right\rvert$. Since $u^{\ast}_{\mu} = \left\lvert u \right\rvert^{\ast}_{\mu}$, the proposition holds also for the function $u$.
\end{myproof4}

\section{Examples}\label{ex}
In this section we will show some examples. We will describe the optimal space $Z'_{m}(\Sigma, \mu)$ when $X(\Sigma, \mu)$ is a Lorentz--Karamata space $L^{p, q, b}(\Sigma, \mu)$ with $p \in \left[1, \frac{D}{m}\right]$.

Now we introduce the theory of Lorentz--Karamata spaces. For proofs and more information see \cite{MR2091115, EKP:00, MR1927106, P2}. The class of Lorentz--Karamata spaces contain a lot of classical function spaces, such as the Lebesgue spaces, Lorentz spaces, Lorentz-Zygmund spaces or Orlicz spaces of the type $L_{\exp}$ or $L^{p}(\log L)^{\alpha}$.

We start with the concept of slowly varying functions. Let $f, g \colon (0, \infty) \to (0, \infty)$ be functions. We say that the function $f$ is \emph{equivalent} to the function $g$ if there exist constants $C_{1} > 0$ and $C_{2} > 0$ such that $C_{2} g(t) \leq f(t) \leq C_{1} g(t)$ for every $t \in (0, \infty)$.

We say that a locally absolutely continuous function $b \colon (0, \infty) \rightarrow (0, \infty)$ is \emph{slowly varying} if for every $\varepsilon > 0$ there exists a nondecreasing function $\varphi_{\varepsilon}$ and a nonincreasing function $\varphi_{- \varepsilon}$ such that $t^{\varepsilon}b(t)$ is equivalent to $\varphi_{\varepsilon}(t)$ on $(0, \infty)$ and that $t^{- \varepsilon}b(t)$ is equivalent to $\varphi_{- \varepsilon}$ on $(0, \infty)$. A positive locally absolutely continuous function on $(0, \infty)$ that is equivalent to a positive constant function is a trivial example of a slowly varying function. Functions of logarithmic type constitute less trivial and very important examples of slowly varying functions. For $k\in\N$, the function $\ell_k\colon (0, \infty) \to (0, \infty)$ defined as
\begin{equation*}
\ell_k(t) = \begin{cases}
1 + |\log t| \quad &\text{if $k = 1$}, \\
1 + \log \ell_{k - 1}(t) \quad &\text{if $k > 1$},
\end{cases}
\end{equation*}
$t\in(0, \infty)$, is slowly varying. More generally, the function $\ell_k^{\mathbb{A}}\colon (0, \infty) \to (0, \infty)$ defined as
\begin{equation*}
\ell_k^{\mathbb{A}}(t) = \begin{cases}
\ell_k^{\alpha_0}(t) \quad &\text{if $t\in(0, 1)$}, \\
\ell_k^{\alpha_\infty}(t) \quad &\text{if $t\in[1, \infty)$},
\end{cases}
\end{equation*}
where $\mathbb {A} = (\alpha_0, \alpha_\infty)\in\R^2$, is slowly varying.

Note that by $b^{-1}$ we will denote the function $\frac{1}{b}$.

Now we finally introduce Lorentz--Karamata spaces. Let $b$ be a slowly varying function. Let $p, q  \in [1, \infty]$. We define the Lorentz--Karamata functionals by
\begin{align*}
\left\lVert f \right\rVert _{p, q, b} = \left\lVert t^{\frac{1}{p}-\frac{1}{q}} b(t) f^{\ast}_{\mu}(t) \right\rVert_{L^{q}(0, \infty)}
\end{align*}
and by
\begin{align*}
\left\lVert f \right\rVert _{(p, q, b)} = \left\lVert t^{\frac{1}{p}-\frac{1}{q}} b(t) f^{\ast \ast}_{\mu}(t) \right\rVert_{L^{q}(0, \infty)}
\end{align*}
for every $f \in \mathcal{M}(\Sigma, \mu)$. The \emph{Lorentz--Karamata spaces} are defined as
\begin{align*}
L^{p, q, b}(\Sigma, \mu) = \left\{f \in \mathcal{M}(\Sigma, \mu); \left\lVert f \right\rVert _{p, q, b} < \infty \right\}
\end{align*}
and as
\begin{align*}
L^{(p, q, b)}(\Sigma, \mu) = \left\{f \in \mathcal{M}(\Sigma, \mu); \left\lVert f \right\rVert _{(p, q, b)} < \infty \right\}.
\end{align*}
If we take $p=q$ and $b\equiv1$, we obtain the Lebesgue spaces. More generally, we obtain the Lorentz spaces $L^{p,q}(\Sigma, \mu)$ and $L^{(p,q)}(\Sigma, \mu)$ by taking $b\equiv 1$. Furthermore, Lorentz--Karamata spaces also include the Lorentz--Zygmund spaces, which were thoroughly studied in \cite{OP:99}. We obtain them by taking slowly varying functions of logarithmic type as above.

Even though we refer to Lorentz--Karamata spaces as spaces, they are not always re\-ar\-range\-ment\--in\-vari\-ant spaces. 

The space $L^{(p, q, b)}(R, \mu)$ is a re\-ar\-range\-ment\--in\-vari\-ant Banach function space if and only if $q\in[1, \infty]$ and one of the following conditions holds:
\begin{enumerate}
\item $p \in (1, \infty)$,
\item $p = 1$ and $\left\lVert t^{-\frac{1}{q}} b(t) \chi_{(1, \infty)}(t) \right\rVert_{L^{q}(0, \infty)} < \infty$,
\item $p = \infty$ and $\left\lVert t^{-\frac{1}{q}} b(t) \chi_{(0, 1)}(t) \right\rVert_{L^{q}(0, \infty)} < \infty$. 
\end{enumerate}

The Lorentz--Karamata functional $\left\lVert \cdot \right\rVert_{p, q, b}$ is equivalent to a~re\-ar\-range\-ment\--in\-vari\-ant Banach function norm if and only if $q\in[1, \infty]$ and one of the following conditions is satisfied:
\begin{enumerate}
\item $p \in (1, \infty)$,
\item $p = q =1$ and $b$ is equivalent to a nonincreasing function on $(0, \infty)$,
\item $p = \infty$ and $\left\lVert t^{-\frac{1}{q}} b(t) \chi_{(0, 1)}(t) \right\rVert_{L^{q}(0, \infty)} < \infty$.
\end{enumerate}
In what follows we will also use the fact that if $p > 1$, then the spaces $L^{p, q, b}(\Sigma, \mu)$ and $L^{(p, q, b)}(\Sigma, \mu)$ are equivalent.

Now we describe the optimal target space for the Lorentz--Karamata space $L^{p, q, b}(\Sigma, \mu)$.
\begin{thm}\label{sccth}
Let $p \in \left [1, \frac{D}{m}\right ]$, let $q \in [1, \infty]$ and let $b$ be a slowly varying function. Assume that $X(\Sigma, \mu) = L^{p, q, b}(\Sigma, \mu)$.
Then we have the following description of the optimal target space $Z'_{m}X(\Sigma, \mu)$.
\begin{enumerate}
\item $Z'_{m}X(\Sigma, \mu) = L^{\frac{Dp}{D-mp}, q, b}(\Sigma, \mu)$ if $p \in \left (1, \frac{D}{m} \right)$ or $p = q = 1$ and $b$ is equivalent to a nonincreasing function on $(0, \infty)$.
\item $Z'_{m}X(\Sigma, \mu) = L^{\infty, q, a}(\Sigma, \mu)$ if $p = \frac{D}{m}$, $q \in (1, \infty]$ and it holds that
\begin{align}\label{asuq}
\left\lVert t^{-\frac{1}{q'}}b^{-1}(t) \right\rVert_{L^{q'}(1,\infty)} < \infty.
\end{align}
Here, the function $a$ is defined by
\begin{align*}
a(t) = \frac{b^{1-q'}(t)}{\int_{t}^{\infty} \tau^{-1}b^{-q'}(\tau) \, d\tau}, \, \, \, \, t \in (0, \infty).
\end{align*} 
\item $Z'_{m}X(\Sigma, \mu) = \Lambda^1(d')(\Sigma, \mu)$ if $p = \frac{D}{m}$, $q = 1$, (\ref{asuq}) is satisfied and $\lim\limits_{t\to 0_+} d(t) = 0$. Here, the function $d$ is defined by 
\begin{align}\label{ddefd}
d(t) = \inf\limits_{\tau \in [t, \infty)} b(\tau), \, \, \, \, t \in (0, \infty).
\end{align}
\item $Z'_{m}X(\Sigma, \mu) = \Lambda^1(d')(\Sigma, \mu) \cap L^\infty(\Sigma, \mu)$ if $p = \frac{D}{m}$, $q = 1$, (\ref{asuq}) is satisfied and $\lim\limits_{t\to 0_+} d(t) > 0$. The function $d$ is again defined by \eqref{ddefd}.
\end{enumerate}
\end{thm}

 The function space $\Lambda^1(d')(\Sigma, \mu)$ is an example of a classical Lorentz space (e.g., see \cite[Chapter~10]{PKJFbook}). It is defined as the collection of all functions $f\in\mathcal M(\Sigma, \mu)$ such that $\|d'(t)f_{\mu}^*(t)\|_{L^1(0, \infty)} < \infty$.

\begin{myproof}
First of all we prove that the space $Z'_{m}X(\Sigma, \mu)$ is the optimal space. Owing to \cite[Section 3.7]{P2} we obtain the fact that the space $\left(L^{p, q, b}(\Sigma, \mu) \right)'$ is equivalent to the space $L^{p', q', b^{-1}}(\Sigma, \mu)$. So, in the view of Theorem~\ref{zopt}, it is enough to prove that 
\begin{align}\label{chcnm}
\left\lVert t^{\frac{m}{D}-1}\chi_{(1, \infty)}(t) \right\rVert_{L^{p', q', b^{-1}}(0, \infty)} = \left\lVert t^{\frac{1}{p'}-\frac{1}{q'}} (t+1)^{\frac{m}{D}-1} b^{-1}(t) \right\rVert_{L^{q'}(0, \infty)} < \infty.
\end{align}
If $p \in [1, \frac{D}{m})$, then (\ref{chcnm}) follows from the basic properties of slowly varying functions (see \cite[Lemma 2.16]{P2}). If $p = \frac{D}{m}$ then (\ref{chcnm}) follows from \eqref{asuq} and from Theorem~\ref{zopt}.

It remains to describe the space $Z'_{m}X(\Sigma, \mu)$. Firstly, assume that $p \in \left (1, \frac{D}{m} \right)$. Owing to \cite[Section 3.5]{P2} the maximal nonincreasing operator $P_{\mu}$ is bounded on the space $L^{p, q, b}(0, \infty)$. Now we show that the norm of the space $Z_{m}X(\Sigma, \mu)$ is equivalent to the functional 
\begin{align}\label{fcl1}
u \mapsto \left\lVert t^{\frac{1}{q}-\frac{1}{p}} b^{-1}(t) \int_{t}^{\infty} u^{\ast \ast}_{\mu} (\tau) \tau^{\frac{m}{D}-1} \, d\tau \right\rVert_{L^{q'}(0, \infty)}, \, \, \, \, u \in \mathcal{M}^{+}(\Sigma, \mu).
\end{align}
By virtue of \cite[Section 3.7]{P2} we can rewrite the functional in the form
\begin{align*}
\begin{split}
& \left\lVert \int_{t}^{\infty}u^{\ast \ast}_{\mu}(\tau) \tau^{\frac{m}{D}-1} \, d\tau \right\rVert_{(L^{p, q, b})'(0, \infty)} \\ & = \sup\limits_{\left\lVert g \right\rVert_{L^{p, q, b}(0, \infty)} \leq 1} \int_{0}^{\infty} g^{\ast}(t) \int_{t}^{\infty}u^{\ast \ast}_{\mu}(\tau) \tau^{\frac{m}{D}-1} \, d\tau dt \\ & = \sup\limits_{\left\lVert g \right\rVert_{L^{p, q, b}(0, \infty)} \leq 1} \int_{0}^{\infty}u^{\ast \ast}_{\mu}(\tau) \tau^{\frac{m}{D}} g^{\ast \ast}(\tau) \, d\tau.
\end{split}
\end{align*}
By virtue of (\ref{cpsp}) we obtain the first inequality
\begin{align*}
\begin{split}
& \left\lVert u \right\rVert_{Z_{m}X(\Sigma, \mu)} = \left\lVert t^{\frac{m}{D}} u^{\ast \ast}_{\mu}(t) \right\rVert_{(L^{p, q, b})'(0, \infty)} \\ & \leq 4 \sup\limits_{\left\lVert g \right\rVert_{L^{p, q, b}(0, \infty)} \leq 1} \int_{0}^{\infty}  g^{\ast \ast}(\tau) \tau^{\frac{m}{D}} u^{\ast \ast}_{\mu}(\tau) \, d\tau,
\end{split}
\end{align*}
while the second inequality follows directly from the boundedness of the operator $P_{\mu}$.
The norm of the associate space $\left( L^{\frac{Dp}{D-mp}, q, b}(\Sigma, \mu) \right)'$ is equivalent to the functional
\begin{align}\label{fcl2}
u \mapsto \left\lVert t^{\frac{1}{q} - \frac{1}{p} + \frac{m}{D}} b^{-1}(t) u^{\ast \ast}_{\mu}(t) \right\rVert_{L^{q'}(0, \infty)},  \, \, \, \, u \in \mathcal{M}^{+}(\Sigma, \mu),
\end{align}
owing to \cite[Section 3.7]{P2}. It remains to prove that the functionals (\ref{fcl1}) and (\ref{fcl2}) are equivalent. Since
\begin{align*}
t^{\frac{m}{D}} = \frac{m}{D \left(1-2^{-\frac{m}{D}}\right)} \int_{\frac{t}{2}}^{t} \tau^{\frac{m}{D}-1} \, d\tau,
\end{align*}
we obtain
\begin{align}\label{tmd}
\begin{split}
& t^{\frac{m}{D}}  u^{\ast \ast}_{\mu}(t) \leq \frac{m}{D \left(1-2^{-\frac{m}{D}}\right)} \int_{\frac{t}{2}}^{t}  u^{\ast \ast}_{\mu}(\tau) \tau^{\frac{m}{D}-1} \, d\tau \\ & \leq \frac{m}{D \left(1-2^{-\frac{m}{D}}\right)} \int_{\frac{t}{2}}^{\infty}  u^{\ast \ast}_{\mu}(\tau) \tau^{\frac{m}{D}-1} \, d\tau
\end{split}
\end{align}
for every $u \in \mathcal{M}^{+}(\Sigma, \mu)$ and for every $t \in (0, \infty)$. By virtue of \eqref{tmd} and \cite[Lemma~2.16]{P2} , we obtain existence of a positive constant $C_{1} $ such that
\begin{align}\label{tmd2}
t^{\frac{1}{q} - \frac{1}{p} + \frac{m}{D}} b^{-1}(t)  u^{\ast \ast}_{\mu}(t) \leq \frac{2^{\frac{1}{q}-\frac{1}{p}}mC_{1}}{D \left(1-2^{-\frac{m}{D}}\right)} \left( \frac{t}{2} \right)^{\frac{1}{q}-\frac{1}{p}} b^{-1}\left(\frac{t}{2} \right) \int_{\frac{t}{2}}^{\infty}  u^{\ast \ast}_{\mu}(\tau) \tau^{\frac{m}{D}-1} \, d\tau
\end{align}
for every $u \in \mathcal{M}^{+}(\Sigma, \mu)$ and for every $t \in (0, \infty)$. Owing to (\ref{dila}) and (\ref{tmd2}) we obtain existence of a positive constant $C_{2}$ such that
\begin{align*}
\left\lVert t^{\frac{1}{q} - \frac{1}{p} + \frac{m}{D}} b^{-1}(t)  u^{\ast \ast}_{\mu} \right\rVert_{L^{q'}(0, \infty)} \leq C_{2} \left\lVert t^{\frac{1}{q}-\frac{1}{p}} b^{-1}(t) \int_{t}^{\infty}  u^{\ast \ast}_{\mu} (\tau) \tau^{\frac{m}{D}-1} \, d\tau \right\rVert_{L^{q'}(0, \infty)}
\end{align*}
for every $u \in \mathcal{M}^{+}(\Sigma, \mu)$. Now we prove the opposite inequality, i.e., we prove that there exists a constant $K_{1} > 0$ such that
\begin{align}\label{nrce}
\begin{split}
& \left\lVert t^{\frac{1}{q}-\frac{1}{p}} b^{-1}(t) \int_{t}^{\infty} u^{\ast \ast}_{\mu} (\tau) \tau^{\frac{m}{D}-1} \, d\tau \right\rVert_{L^{q'}(0, \infty)} \\ & \leq K_{1} \left\lVert t^{\frac{1}{q} - \frac{1}{p} + \frac{m}{D}} b^{-1}(t) u^{\ast \ast}_{\mu}(t) \right\rVert_{L^{q'}(0, \infty)},
\end{split}
\end{align}
for every $u \in \mathcal{M}^{+}(\Sigma, \mu)$. By virtue of the weighted Hardy inequality (\cite[Theorem 2]{MHP}), the inequality (\ref{nrce}) is valid if and only if
\begin{align*}
\sup\limits_{t \in (0, \infty)} \left\lVert \tau^{\frac{1}{q}-\frac{1}{p}} b^{-1}(\tau) \right\rVert_{L^{q'}(0, t)} \left\lVert \tau^{\frac{1}{p} - \frac{1}{q} - 1} b(\tau) \right\rVert_{L^{q}(t, \infty)} < \infty, 
\end{align*}
which is not hard to verify (see~\cite[Lemma 2.16]{P2}).

It remains to derive the description of the space $Z'_{m}X(\Sigma, \mu)$ in case of $p = q = 1$ and $b$ is equivalent to a nonincreasing function on $(0, \infty)$ or $ p = \frac{D}{m}$ and $q \in [1, \infty]$. Thanks to \cite[Section 3.5 and Section 3.7]{P2} it is sufficient to prove that $Z_{m}X(\Sigma, \mu)$ with the norm satisfying
\begin{align*}
\left\lVert u \right\rVert_{Z_{m}X(\Sigma, \mu)} = \left\lVert t^{\frac{m}{D}} u^{\ast \ast}_{\mu}(t) \right\rVert_{L^{p', q', b^{-1}}(0, \infty)}
\end{align*}
for each $u \in \mathcal{M}^{+}(\Sigma, \mu)$ is equivalent to $L^{(r, q', b^{-1})}(\Sigma, \mu)$, where $r = 1$ if $p = \frac{D}{m}$ and $r = \frac{D}{m}$ if $ p = 1$. Note that in case of $p = \frac{D}{m}, q = 1$, we also exploit the fact that $d^{-1}$ is locally absolutely continuous on $(0, \infty)$. Now, for every $u \in \mathcal{M}^{+}(\Sigma, \mu)$ we have
\begin{align*}
\begin{split}
& \left\lVert u\right\rVert_{Z_{m}X(\Sigma, \mu)} = \left\lVert t^{\frac{m}{D}} u^{\ast \ast}_{\mu}(t) \right\rVert_{L^{p', q', b^{-1}}(0, \infty)} \leq \left\lVert \sup\limits_{\tau \in [t, \infty)} \tau^{\frac{m}{D}} u^{\ast \ast}_{\mu}(\tau) \right\rVert_{L^{p', q', b^{-1}}(0, \infty)} \notag \\ & = \left\lVert t^{\frac{1}{p'} - \frac{1}{q'}} b^{-1}(t) \sup\limits_{\tau \in [t, \infty)} \tau^{\frac{m}{D}} u^{\ast \ast}_{\mu}(\tau) \right\rVert_{L^{q'}(0, \infty)} \\ & \leq K_{2} \left\lVert t^{\frac{1}{r} - \frac{1}{q'}} b^{-1}(t) u^{\ast \ast}_{\mu}(t) \right\rVert_{L^{q'}(0, \infty)} \notag = \left\lVert u \right\rVert_{L^{(r, q', b^{-1})}(\Sigma, \mu)},
\end{split} 
\end{align*}
where $K_{2}$ is a positive constant that does not depend on $u$. If $p = q = 1$, then the last inequality follows from the fact that $b^{-1}$ is equivalent to a nonincreasing function on $(0, \infty)$. If $p = \frac{D}{m}, q  = 1$, then the last inequality follows from the fact that $t^{1 - \frac{m}{D}}b^{-1}(t)$ is equivalent to a nonincreasing function on $(0, \infty)$. If $q > 1$, then we use \cite[Theorem 3.2]{OPG} to prove the last inequality. Owing to this theorem the last inequality is valid if and only if for every $t \in (0, \infty)$ it is true that
\begin{align*}
t^{\frac{m}{D}} \left\lVert \tau^{\frac{D-m}{D} - \frac{1}{q'}} b^{-1}(\tau) \right\rVert_{L^{q'}(0, t)} \leq K \left\lVert \tau^{1 - \frac{1}{q'}} b^{-1}(\tau) \right\rVert_{L^{q'}(0, t)},
\end{align*}
where $K > 0$ is a constant that does not depend on $t$, which is satisfied thanks to \cite[Lemma 2.16]{P2}. So, we have proved the first inequality. For every $u \in \mathcal{M}^{+}(\Sigma, \mu)$ we have now
\begin{align*}
\begin{split}
& \left\lVert u \right\rVert_{L^{(r, q', b^{-1})}(\Sigma, \mu)} = \left\lVert t^{\frac{1}{r} - \frac{1}{q'}} b^{-1}(t) u^{\ast \ast}_{\mu}(t) \right\rVert_{L^{q'}(0, \infty)} \\ & \leq \left\lVert t^{\frac{1}{r} - \frac{m}{D} - \frac{1}{q'}} b^{-1}(t) \sup\limits_{\tau \in [t, \infty)} \tau^{\frac{m}{D}} u^{\ast \ast}_{\mu}(\tau) \right\rVert_{L^{q'}(0, \infty)} \\ & = \left\lVert \sup\limits_{\tau \in [t, \infty)} \tau^{\frac{m}{D}} u^{\ast \ast}_{\mu}(\tau) \right\rVert_{L^{p', q', b^{-1}}(0, \infty)} \leq C \left\lVert t^{\frac{m}{D}} u^{\ast \ast}_{\mu}(t) \right\rVert_{L^{p', q', b^{-1}}(0, \infty)} \\ & = C \left\lVert u\right\rVert_{Z_{m}X(\Sigma, \mu)},
\end{split}
\end{align*}
where $C$ is a positive constant that does not depend on $u$. To prove the second inequality, we exploit \cite[Lemma 4.10]{EMMP} since the function $t \mapsto t^{\frac{m}{D}}, t \in [0, \infty),$ is quasiconcave.
\end{myproof}

\begin{rem}
If the function $b$ is in addition nonincreasing on $(0, 1]$ and constant on $[1, \infty)$, then we obtain the fact that $d$ is constant on $(0, \infty)$, so $d' \equiv 0$ on $(0, \infty)$. It means that the optimal space $Z'_{m}X(\Sigma, \mu)$ is equivalent to $L^{\infty}(\Sigma, \mu)$. For instance, if $d \equiv 1$ on $(0, \infty)$, then we obtain the Lorentz space $L^{\frac{D}{m}, 1}(\Sigma, \mu)$ (cf.~\cite{S:81}). 
\end{rem}

Finally, we describe the optimal domain space in the view of Theorem~\ref{zopt22} for a given Lorentz--Karamata space.

\begin{thm}\label{sccth1}
Let $p \in \left [\frac{D}{D-m}, \infty \right ]$, let $q \in [1, \infty]$ and let $b$ be a slowly varying function. Assume that $Y(\Sigma, \mu) = L^{p, q, b}(\Sigma, \mu)$. Then we have the following description of the optimal domain space $U_{m}Y(\Sigma, \mu)$.
\begin{enumerate}
\item $U_{m}Y(\Sigma, \mu) = L^{\frac{Dp}{D + mp}, q, b}(\Sigma, \mu)$ if $p \in \left (\frac{D}{D-m}, \infty \right)$.
\item $U_{m}Y(\Sigma, \mu) = L^{1,1,b}(\Sigma, \mu)$ if $p = \frac{D}{D-m}$, $q = 1$ and $b$ is equivalent to a nonincreasing function on $(0, \infty)$.
\item\label{drei} $U_{m}Y(\Sigma, \mu)$ is given by
\begin{align*}
 \left\lVert f \right\rVert_{U_{m}Y(\Sigma, \mu)} = \left\lVert t^{-\frac{1}{q}} b(t) \int_{t}^{\infty} f^{\ast}_{\mu}(\tau) \tau^{\frac{m}{D}-1} \, d\tau \right\rVert_{L^{q}(0, \infty)}
\end{align*}
if $p = \infty$ and
\begin{align*}
\left\lVert t^{-\frac{1}{q}}b(t) \right\rVert_{L^{q}(0,1)} < \infty.
\end{align*}
\end{enumerate}
\end{thm}

\begin{myproof}
For every $f \in \mathcal{M}(0, \infty)$, we define the operator $T_{\frac{m}{D}}(f)$ by
\begin{align*}
T_{\frac{m}{D}}(f)(t) = t^{-\frac{m}{D}} \sup\limits_{\tau \in [t, \infty)} \tau^{\frac{m}{D}} f^{\ast}_{\mu}(\tau), \,\,\,\, t \in (0, \infty).
\end{align*}
Owing to \cite[Theorem 4.7 and Remark 4.8]{EMMP}, we obtain the fact that the operator $T_{\frac{m}{D}} \colon Y'(0, \infty) \to Y'(0, \infty)$ is bounded if and only if the norm of the space $W_{m}Y(\Sigma, \mu)$ is equivalent to 
\begin{align}\label{eqopt}
\left\lVert \int_{t}^{\infty} f^{\ast}_{\mu}(\tau) \tau^{\frac{m}{D}-1} \, d\tau \right\rVert_{Y(0, \infty)}.
\end{align}

So, we prove the fact that the operator $T_{\frac{m}{D}}$ is bounded on the space $\left(L^{p, q, b}(0, \infty)\right)'$. Firstly, assume that $p < \infty$ or $q < \infty$. We use \cite[Section 3.7]{P2} to describe the space $\left(L^{p, q, b}(0, \infty)\right)'$. If $p \in \left [\frac{D}{D-m}, \infty \right) $, then we obtain the fact that $\left(L^{p, q, b}(\Sigma, \mu) \right)'$ is equivalent to $ L^{(p', q', b^{-1})}(\Sigma, \mu)$. If $p = \infty$ and $q \in [1, \infty)$, then $\left(L^{\infty, q, b}(\Sigma, \mu) \right)'$ is equivalent to $L^{(1, q', a)}(\Sigma, \mu)$, where $a$ is the slowly varying function given by
\begin{align*}
a(t) = \left( \int\limits_{0}^{t} \tau^{-1} b^{q}(\tau) \, d\tau \right)^{-1} b^{q-1}(t), \,\,\,\, t \in (0, \infty).
\end{align*}
In what follows we will write $c$ instead of $b^{-1}$ and $a$ since the estimates are the same in both cases. Then, for every $f \in \mathcal{M}(0, \infty)$, we obtain
\begin{align}\label{lasst}
\begin{split}
& \left\lVert T_{\frac{m}{D}}(f) \right\rVert_{L^{(p', q', c)}(0, \infty)} = \left\lVert t^{\frac{1}{p'} - \frac{1}{q'}} c(t) \left(T_{\frac{m}{D}}(f)\right)^{\ast \ast}_{\mu}(t) \right\rVert_{L^{q'}(0, \infty)} \\ & \leq C_{2} \left\lVert t^{\frac{1}{p'} - \frac{1}{q'}} c(t) T_{\frac{m}{D}}(f^{\ast \ast}_{\mu})(t) \right\rVert_{L^{q'}(0, \infty)} \\ & = C_{2} \left\lVert t^{\frac{1}{p'} - \frac{1}{q'}} c(t) t^{-\frac{m}{D}} \sup\limits_{\tau \in [t, \infty)} \tau^{\frac{m}{D}} f^{\ast \ast}_{\mu}(\tau) \right\rVert_{L^{q'}(0, \infty)} \\ & \leq C_{3} \left\lVert t^{\frac{1}{p'} - \frac{1}{q'}} c(t) f^{\ast \ast}_{\mu}(t) \right\rVert_{L^{q'}(0, \infty)} = C_{3} \left\lVert f \right\rVert_{L^{(p', q', c)}(0, \infty)}.
\end{split}
\end{align}
The first inequality is true thanks to \cite[Lemma 4.1]{MOS}. Now we prove the second inequality. If $q > 1$, then we use \cite[Theorem 3.2]{OPG} to prove it. Owing to this theorem the second inequality is valid if and only if for every $t \in (0, \infty)$ it is true that
\begin{align}\label{hrdx}
t^{\frac{m}{D}} \left\lVert \tau^{\frac{1}{p'} - \frac{1}{q'} - \frac{m}{D}} c(\tau) \right\rVert_{L^{q'}(0, t)} \leq K \left\lVert \tau^{\frac{1}{p'} - \frac{1}{q'}} c(\tau) \right\rVert_{L^{q'}(0, t)},
\end{align}
where $K > 0$ is a constant that does not depend on $t$, which is not hard to verify (see~\cite[Lemma 2.16]{P2}). If $q = 1$, then $q' = \infty$, so the second inequality in \eqref{lasst} is of the form
\begin{align*}
\begin{split}
& C_{2} \sup\limits_{t \in (0, \infty)} t^{\frac{1}{p'} - \frac{m}{D}} c(t) \sup\limits_{\tau \in [t, \infty)} \tau^{\frac{m}{D}} f^{\ast \ast}_{\mu}(\tau) = C_{2} \sup\limits_{\tau \in (0, \infty)} \tau^{\frac{m}{D}} f^{\ast \ast}_{\mu}(\tau) \sup\limits_{t \in (0, \tau]}  t^{\frac{1}{p'} - \frac{m}{D}} c(t) \\ & \leq C_{3} \sup\limits_{t \in (0, \infty)} t^{\frac{1}{p'}} c(t) f^{\ast \ast}_{\mu}(t) = C_{3} \left\lVert f \right\rVert_{L^{(p', \infty, c)}(0, \infty)}.
\end{split}
\end{align*}
Note that we used the fact that $b^{-1}$ is equivalent to a nondecreasing function on $(0, \infty)$ if $p = \frac{D}{D-m}$. It means that \eqref{lasst} is true. So, the operator $T_{\frac{m}{D}}$ is bounded on the space $\left(L^{p, q, b}(0, \infty)\right)'$. If $p = q = \infty$, then we use \cite[Proposition 3.3 and Theorem 3.31]{P2} to obtain the fact that $\left(L^{\infty, \infty, b}(\Sigma, \mu) \right)'$ is equivalent to $L^{1, 1, h^{-1}}(\Sigma, \mu)$, where $h$ is the slowly varying function given by $h(t) = \sup\limits_{\tau \in (0, t)} b(\tau)$, $t \in (0, \infty)$. The boundedness of the operator $T_{\frac{m}{D}}$ in this case can be proved similarly as in the preceding cases hence we omit the proof. Now it follows from the proof of \cite[Theorem 4.1]{EMMP} that the condition \eqref{jnvm22} is satisfied. It means that the norm of the optimal space $U_{m}Y(\Sigma, \mu)$ given by Theorem~\ref{zopt22} is equivalent to \eqref{eqopt}. So, we have proved the assertion \eqref{drei} of the theorem. 

It remains to prove that \eqref{eqopt} is equivalent to the norm of the space $U_{m}Y(\Sigma, \mu)$ in the first two cases. For each $f \in \mathcal{M}(0, \infty)$ we obtain
\begin{align*}
\begin{split}
& \left\lVert \int_{t}^{\infty} f^{\ast}_{\mu}(\tau) \tau^{\frac{m}{D}-1} \, d\tau \right\rVert_{L^{p, q, b}(0, \infty)} = \left\lVert t^{\frac1{p} - \frac1{q}} b(t) \int_{t}^{\infty} f^{\ast}_{\mu}(\tau) \tau^{\frac{m}{D}-1} \, d\tau \right\rVert_{L^{p, q, b}(0, \infty)} \\
&\geq \left\lVert t^{\frac1{p} - \frac1{q}} b(t) \int_{t}^{2t} f^{\ast}_{\mu}(\tau) \tau^{\frac{m}{D}-1} \, d\tau \right\rVert_{L^{p, q, b}(0, \infty)} \geq C_{4} \left\lVert t^{\frac{1}{p} - \frac{1}{q} + \frac{m}{D}} b(t) f^{\ast}_{\mu}(2t) \right\rVert_{L^{q}(0, \infty)} \\
&\geq C_{5} \left\lVert t^{\frac{1}{p} - \frac{1}{q} + \frac{m}{D}} b(t) f^{\ast}_{\mu}(t) \right\rVert_{L^{q}(0, \infty)} = C_{5}\left\lVert f \right\rVert_{L^{\frac{Dp}{D + mp}, q, b}(0, \infty)}.
\end{split}
\end{align*}
In the third inequality, we used the boundedness of the dilation operator (see \eqref{dila}). On the other hand, for every $f \in \mathcal{M}(0, \infty)$ we have
\begin{align*}
\begin{split}
& \left\lVert \int_{t}^{\infty} f^{\ast}_{\mu}(\tau) \tau^{\frac{m}{D}-1} \, d\tau \right\rVert_{L^{p, q, b}(0, \infty)} = \left\lVert t^{\frac{1}{p} - \frac{1}{q}} b(t)  \int_{t}^{\infty} f^{\ast}_{\mu}(\tau) \tau^{\frac{m}{D}-1} \, d\tau \right\rVert_{L^{q}(0, \infty)} \\ & \leq C_{6} \left\lVert t^{\frac{1}{p} - \frac{1}{q} + \frac{m}{D}} b(t)  f^{\ast}_{\mu}(t)\right\rVert_{L^{q}(0, \infty)} = C_{6} \left\lVert f \right\rVert_{L^{\frac{Dp}{D + mp}, q, b}(0, \infty)}.
\end{split}
\end{align*}
By virtue of the weighted Hardy inequality (\cite[Theorem 2]{MHP}), the preceding inequality is valid if and only if
\begin{align*}
\sup\limits_{t \in (0, \infty)} \left\lVert \tau^{\frac{1}{p}-\frac{1}{q}} b(\tau) \right\rVert_{L^{q}(0, t)} \left\lVert \tau^{\frac{1}{q} - \frac{1}{p} - 1} b^{-1}(\tau) \right\rVert_{L^{q'}(t, \infty)} < \infty, 
\end{align*}
which is not hard to verify (see~\cite[Lemma 2.16]{P2}).
\end{myproof}

\begin{rem}
It is a natural question what we can say about the optimal domain space for the Lorentz--Karamata space $Y(\Sigma, \mu)= L^{p, q, b}(\Sigma, \mu)$ in the case $p = \frac{D}{D-m}, q \in (1, \infty]$. In this case we have 
\begin{align}\label{nrmr}
\frac{t^{1-\frac{m}{D}}}{\varphi_{Y}(t)} = \frac{t^{1-\frac{m}{D}}}{\left\lVert \chi_{(0, t)} \right\rVert_{L^{\frac{D}{D-m}, q, b}(0, \infty)}} = \frac{t^{1-\frac{m}{D}}}{\left\lVert \tau^{1 - \frac{m}{D} - \frac{1}{q}} b(\tau)  \right\rVert_{L^{ q }(0, t)}}.
\end{align}
It can be proved by \cite[Lemma 2.16]{P2} and elementary computations that the right-hand side of \eqref{nrmr} is equivalent to $b^{-1}$ on $(0, \infty)$. It means that the condition \eqref{jnvm22} is satisfied if and only if $b$ is equivalent to a nonincreasing function on $(1, \infty)$. By Theorem~\ref{zopt22} we conclude that the optimal domain space for $Y(\Sigma, \mu) = L^{\frac{D}{D-m}, q, b}(\Sigma,\mu)$ exists if and only if $b$ is equivalent to a nonincreasing function on $(1, \infty)$. Furthermore, if this is not the case, then there is no domain space $X$ with which \eqref{iqsy22} is valid.

On the other hand, it can be proved by \cite[Theorem~3.2]{OPG} that the operator $T_{\frac{m}{D}}$ is not bounded on the space $\left(L^{\frac{D}{D-m}, q, b}\right)'(0, \infty)$, which is equivalent to $L^{\frac{D}{m}, q', b^{-1}}(0, \infty)$ owing to \cite[Section~3.7]{P2}. In particular, notice that the condition \eqref{hrdx} is not satisfied, which is not hard to verify thanks to \cite[Lemma~2.16]{P2}). The unboundedness of $T_{\frac{m}{D}}$ on $\left(L^{\frac{D}{D-m}, q, b}\right)'(0, \infty)$ means that we cannot simplify the norm of the optimal space $U_{m}Y(\Sigma, \mu)$ (\cite[Theorem~4.7 and Remark~4.8]{EMMP}).
\end{rem}
\bibliography{bibliography}
\end{document}